\documentclass[12pt,a4paper]{amsart}
\usepackage{amsmath, amssymb, amsthm,
amsfonts, amsthm, bm, eucal, graphicx}
\makeatletter
\@namedef{subjclassname@2020}{\textup{2020} Mathematics Subject Classification}
\makeatother
\usepackage[all]{xy}
\usepackage{graphics}
\usepackage{epsfig}
\usepackage{graphicx}
\theoremstyle{plain}
\usepackage{amssymb}
\usepackage[colorlinks]{hyperref}
\usepackage{enumerate}

\advance\hoffset-20mm \advance\textwidth40mm

\newtheorem{thm}{Theorem}[section]
\newtheorem{cor}[thm]{Corollary}
\newtheorem{lem}[thm]{Lemma}
\newtheorem{prop}[thm]{Proposition}

\theoremstyle{definition}
\newtheorem{defi}[thm]{Definition}
\newtheorem{defis}[thm]{Definitions}
\newtheorem{conj}[thm]{Problem}
\newtheorem{conv}[thm]{Convention}
\newtheorem{nota}[thm]{Notation}
\newtheorem{rem}[thm]{Remark}
\newtheorem{rems}[thm]{Remarks}
\newtheorem{exa}[thm]{Example}
\newtheorem{exas}[thm]{Examples}
\newtheorem{sit}[thm]{}

\newcommand{\brem}{\begin{rem}}
\newcommand{\brems}{\begin{rems}}
\newcommand{\erem}{\end{rem}}
\newcommand{\erems}{\end{rems}}
\newcommand{\bexa}{\begin{exa}}
\newcommand{\bexas}{\begin{exas}}
\newcommand{\eexa}{\end{exa}}
\newcommand{\eexas}{\end{exas}}
\newcommand{\bdefi}{\begin{defi}}
\newcommand{\edefi}{\end{defi}}
\newcommand{\bdefis}{\begin{defis}}
\newcommand{\edefis}{\end{defis}}
\newcommand{\bcor}{\begin{cor}}
\newcommand{\ecor}{\end{cor}}
\newcommand{\blem}{\begin{lem}}
\newcommand{\elem}{\end{lem}}
\newcommand{\bconv}{\begin{conv}}
\newcommand{\econv}{\end{conv}}
\newcommand{\bconj}{\begin{conj}}
\newcommand{\econj}{\end{conj}}
\newcommand{\bprop}{\begin{prop}}
\newcommand{\eprop}{\end{prop}}
\newcommand{\bthm}{\begin{thm}}
\newcommand{\ethm}{\end{thm}}
\newcommand{\bnota}{\begin{nota}}
\newcommand{\enota}{\end{nota}}
\newcommand{\bsit}{\begin{sit}}
\newcommand{\esit}{\end{sit}}
\newcommand{\be}{\begin{equation}}
\newcommand{\ee}{\end{equation}}
\newcommand{\bproof}{\begin{proof}}
\newcommand{\eproof}{\end{proof}}
\def\ba{\begin{array}}
\def\ea{\end{array}}
\def\bea{\begin{eqnarray}}
\def\eea{\end{eqnarray}}

\def\bnum{\begin{enumerate}}
\def\enum{\end{enumerate}}
\newcommand{\la}{\label}

\newtheorem*{theo*}{Theorem}

\theoremstyle{definition}

\newtheorem*{definition*}{Definition}

\def\fG{{\mathfrak G}}

\def\cA{{\mathcal A}}
\def\cB{{\mathcal B}}
\def\kk{{\mathcal C}}

\def\cG{{\mathcal G}}

\def\cN{{\mathcal N}}

\def\cR{{\mathcal R}}

\def\cV{{\mathcal V}}

\def\AA{{\mathbb A}}

\def\NN{{\mathbb N}}
\def\ZZ{{\mathbb Z}}
\def\QQ{{\mathbb Q}}
\def\TT{{\mathbb T}}
\def\kk{{\Bbbk}}
\def\CC{{\mathbb C}}
\def\G{{\mathbb G}}
\def\Ga{{\mathbb G_{\mathrm a}}}
\def\Gm{{\mathbb G_{\mathrm m}}}

\def\RR{{\mathbb R}}
\def\PP{{\mathbb P}}

\def\End{\mathop{\rm End}}
\def\Aut{\mathop{\rm Aut}}
\def\Bir{\mathop{\rm Bir}}

\def\GL{\mathop{\rm GL}}
\def\PGL{\mathop{\rm PGL}}
\def\SL{\mathop{\rm SL}}
\def\Jonq{\mathop{\rm Jonq}}
\def\JONQ{\mathop{\rm JONQ}}

\def\Stab{\mathop{\rm Stab}}
\def\Norm{\mathop{\rm Norm}}
\def\Cent{\mathop{\rm Cent}}

\def\Spec{\mathop{\rm Spec}}
\def\Trans{\mathop{\rm Trans}}
\def\dist{\mathop{\rm dist}}

\def\deg{\mathop{\rm deg}}

\def\ord{\mathop{\rm ord}}

\def\Aff{\mathop{\rm Aff}}

\def\Ax{\mathop{\rm Ax}}
\def\vert{\mathop{\rm vert}}

\def\ll1{l_{\lambda}^{-1}(1)}
\def\lm1{l_{\mu}^{-1}(1)}

\def\ba{\begin{array}}
\def\ea{\end{array}}
\def\bea{\begin{eqnarray}}
\def\eea{\end{eqnarray}}

\begin{document}
\sloppy
\title[Borel subgroups]
{Borel subgroups of the automorphism 
groups of affine toric surfaces}

\author%
{Ivan\ Arzhantsev and Mikhail\ Zaidenberg}
\address{Faculty of Computer Science, HSE University, 
Pokrovsky Boulevard 11, Moscow,
109028 Russia}
\email{arjantsev@hse.ru}
\address{Univ. Grenoble Alpes, CNRS, IF, 
38000 Grenoble, France}
\email{mikhail.zaidenberg@univ-grenoble-alpes.fr}
%
%
\begin{abstract} 
In \cite{AZ13}
we described the automorphism groups
of the cyclic quotients
of the affine plane $\AA^2$.
In this article, we study
the Borel subgroups and,
more generally, the maximal solvable subgroups
of these ind-groups.
We show that the cyclic quotients
of $\AA^2$ are divided into two species.
In one of them, the Borel subgroups
form a single conjugacy class,
while in the other, there are
two conjugacy classes 
of Borel subgroups.
The proofs explore the 
Bass-Serre-Tits theory
of groups acting on trees.
\end{abstract}
\thanks{The work of the first author is supported 
by the grant RSF 25-11-0030.}
\subjclass[2020]{Primary 14J50, 20E08; \ 
Secondary 14M25, 14R10}
\keywords{Affine toric surface, 
automorphism group, Bass--Serre tree, Borel subgroup, maximal solvable subgroup}
\dedicatory{To Leonid Makar-Limanov on occasion of his 80th birthday}
\maketitle
{\footnotesize \tableofcontents}

\section{Introduction}
Let $\kk$ be an algebraically 
closed field  of characteristic zero 
and $X$ an affine 
algebraic variety  over $\kk$.
By definition, a Borel subgroup 
of the automorphism 
group $\Aut(X)$
is a maximal solvable 
connected ind-subgroup 
of $\Aut(X)$.
The classical Borel theorem 
(\cite[Th\'eor\`eme~16.5]{Bor56}) says 
that all Borel subgroups 
of an affine algebraic group are conjugate.
A similar statement holds
for the automorphism group of the
affine plane $\AA^2$ over $\kk$,
see \cite{BEE16} and \cite{FP18}, 
or Theorem \ref{thm:Borel-A2} below. 
Note that $\Aut(\AA^n)$ for $n\ge 3$ 
contains Borel subgroups that are 
not conjugate to the standard Borel subgroup
of triangular automorphisms, 
see \cite[Corollary~1.5]{RUvS24}. 
The same holds for the subgroup 
$\Aut_z(\AA^3)\subset\Aut(\AA^3)$ 
consisting of all automorphisms
of  $\AA^3=\Spec\kk[x,y,z]$ 
fixing the coordinate $z$, 
see \cite[Corollary~3.22]{FP18}.
The Cremona group ${\rm Bir}(\PP^2)$
contains families of 
unbounded dimensions
of pairwise non-conjugate Borel subgroups,  see
\cite[Theorem~1.3]{FH23}. 
Confirming \cite[Conjecture~1.2]{FH23}, 
it has been shown in  \cite[Corollary~1.4]{RUvS24} 
that ${\rm Bir}(\PP^n)$
contains
non-conjugate  
Borel subgroups for every $n\ge 2$.

Let $\zeta\in\mathbb{G}_{\mathrm m}
=\mathbb{G}_{\mathrm m}(\kk)$ 
be a primitive $d$th root of unity.
Our study focuses mainly on the
cyclic quotients 
$X_{d,e}=\AA^2/\mu_d$,
where  
$\mu_d=\langle\zeta\rangle\simeq\ZZ/d\ZZ$ 
acts on the affine plane via
$(x,y)\mapsto (\zeta^e x, \zeta y)$. 
Any normal affine toric surface $X$
different from the affine plane and having
no non-constant invertible regular function
is isomorphic to one of the $X_{d,e}$, 
see \cite[Section~2.6]{Ful93}. 

In this article,
we describe the
Borel subgroups of $\Aut(X_{d,e})$. 
The crucial fact used in the proofs is 
an analogue of the 
Jung--van der Kulk--Nagata
theorem, see \cite[Theorem~4.2]{AZ13}. 
It asserts that
the automorphism groups $\Aut(X_{d,e})$
are amalgams. 

A Borel--Chevalley theorem states that 
every Borel subgroup of 
a connected affine
algebraic group $G$ 
is maximal among the solvable 
subgroups of $G$, 
see \cite[Proposition~18]{Bor56} and 
\cite[Corollaire 2]{Che57}. 
Our main result
is analogues to the Borel and 
Borel--Chevalley theorems. 
Informally, it states the following; 
see Section \ref{sec:toric} 
for the terminology and 
Theorem \ref{thm:Borel} 
for the precise statement.
\bthm\la{main}  $\,$
\bnum
\item[$(a)$]
Let $e^2\equiv 1\mod d$. 
Then $\Aut(X_{d,e})$ has a unique 
 conjugacy
 class of Borel subgroups.
 This class contains a 
 triangular subgroup $\cB$.
Every solvable connected subgroup 
 of $\Aut(X_{d,e})$ is triangulable, i.e., is 
 conjugate to a subgroup of $\cB$. 
 \item[$(b)$]
Let $e^2\not\equiv 1\mod d$.
Then $\Aut(X_{d,e})$ 
has 
exactly two conjugacy classes of Borel subgroups
represented by the triangular 
subgroups $\cB^+$ and $\cB^-$
that are factors of the amalgam 
$\Aut(X_{d,e})=\cB^+*_{\TT}\cB^-$, where 
$\cB^+\cap \cB^-=
\TT\simeq\mathbb{G}^2_{\mathrm m}$ 
is a maximal torus.
Every solvable connected subgroup 
 of $\Aut(X_{d,e})$ is conjugate to a subgroup
of $\cB^+$ or $\cB^-$.
 \item[$(c)$]
Every Borel subgroup of 
$\Aut(X_{d,e})$ is
maximal among the solvable subgroups 
of $\Aut(X_{d,e})$. 
\enum
 \ethm
 Besides the surfaces $\AA^2$ and $X_{d,e}$, 
 there are, up to isomorphism, only two more
 normal affine toric surfaces, namely, $(\AA^1_*)^2$
and $\AA^1\times\AA^1_*$, where
$\AA^1_*=\AA^1\setminus\{0\}$,
see for example \cite[Section~3.2]{KPZ17}.  
 These toric surfaces carry non-constant
 invertible regular functions. 

The automorphism group 
$G=\Aut((\AA^1_*)^2)\simeq
\TT\rtimes\GL(2,\ZZ)$ is an extension of the
2-torus 
$\TT\simeq\mathbb{G}^2_{\mathrm m}$. 
Its unique
Borel subgroup $\TT$ coincides 
with the neutral component 
$G^0$. 

The group $G=\Aut(\AA^1\times\AA^1_*)$
is solvable and also has a unique Borel subgroup 
\[G^0=\{(x,y)\mapsto 
(ax+p(y),by)\,|\,p(y)\in\kk[y,y^{-1}],\,\,\, a,b\in 
\mathbb{G}_{\mathrm m}\},\]
see \cite[Section~3.5]{KPZ17}.

If $X=\AA^2/G$ is the quotient 
of the affine plane by a small finite 
non-abelian subgroup 
$G\subset\GL(2,\kk)$,
then $\Aut(X)$ is a finite extension of the 
1-torus $Z\simeq\Gm$ 
of scalar matrices, 
see \cite[Theorem~5.3]{AZ13} 
(cf. also \cite{Kol24}). 
Thus, the unique Borel subgroup of $\Aut(X)$
coincides with $Z$. 

The following theorem is analogues to 
\cite[Theorem~4]{BEE16}, which is itself analogues to
Steinberg's characterization of Borel subgroups in 
connected affine algebraic groups 
(see \cite[Section~2.1(b)]{Ste74}).
It this theorem,
$\Aut(X_{d,e})$ is considered 
as an abstract group,
without taking into account its 
ind-structure and topology. 
\bthm\la{thm:main-2}
A non-abelian maximal solvable subgroup 
$H$ of $\Aut(X_{d,e})$
is a Borel subgroup 
if and only if it contains no proper 
 subgroup of finite index.
\ethm
Unlike in the case of affine 
algebraic groups, 
the assumption ``non-abelian'' 
in Theorem \ref{thm:main-2}
is essential, as first noted in  \cite{BEE16}
in a similar context. 
For example, 
the automorphism group
$\Aut(\AA^2)$ has abelian maximal 
solvable subgroups  that contain 
no finite-index subgoups 
and yet are not Borel subgroups, 
cf. Remark \ref{rem:Mold-BEE}.3
and Lemma \ref{exa:Wright-BEE}. 

The number of conjugacy classes of 
maximal solvable subgroups in an 
algebraic group defined
over an algebraically closed field
is finite, see \cite[Theorem~A]{Pla69}. 
It is unknown whether an analogue 
of this result 
is true for $\Aut(\AA^n)$, $n\ge 2$. 
As we already mentioned, 
this fails for the Cremona group 
$\Bir(\PP^2)$.

The content of the article is as follows. 
We use
the technique
explored in Lamy's 
proof of Tits' alternative
for $\Aut(\AA^2)$, see \cite{Lam01}. 
It  is based on the 
Bass--Serre--Tits theory 
of groups acting on trees. 
In Section \ref{sec:prelim}, 
we present results from 
the Bass--Serre--Tits theory
that we use in the next sections.
They  mainly concern
the classification of  subgroups 
of an amalgam 
acting on its Bass--Serre tree, 
see Proposition \ref{thm:Tits} 
and Remark \ref{rem:classification}. 
Broadly speaking, all subgroups are divided 
into three main classes: 
elliptic, parabolic, and loxodromic subgroups, 
see Definitions \ref{def:sbgrps}.

In Section \ref{sec:maximal}
we describe conjugacy classes of subgroups 
of an amalgam that are
maximal among the solvable subgroups,
see Propositions \ref{thm:maximal} and
\ref{prop:exas-max-loxo} and 
Corollaries
\ref{cor:lox-solv1}--\ref{cor:lox-solv2}. 
Informally, we can summarize 
these results as follows 
(see Section \ref{sec:prelim} 
for terminology). 
\bthm\la{mthm-3}
Consider an amalgam $G=A*_C B$, 
where $C=A\cap B$ is a proper 
common subgroup 
of $A$ and $B$. Let $T$ be 
the Bass-Serre tree  of $G$.
\bnum
\item[$(a)$] If $A$ is solvable 
and $|A:C| \ge 3$, 
then $A$ is a maximal 
solvable subgroup 
of $G$. 
\item[$(b)$] Suppose that $C$
 is solvable, and let $S$
be a maximal solvable subgroup of $G$. 
Then one of the following statements is true:
\begin{itemize}
\item $S$ is elliptic, that is, 
conjugate to a subgroup 
of $A$ or $B$;
\item $S$ is parabolic and
there exists a Freudenthal 
end $b\in\partial T$ 
such that $S=\Stab_G(b)$;
\item $S$ is loxodromic 
and there exists either
a Freudenthal end $b\in \partial T$ 
such that $S=\Stab_G(b)$,
or a pair 
of distinct  Freudenthal 
ends $b,b'\in\partial T$
such that $S=\Stab_G(\{b,b'\})$.
\end{itemize}
\enum
\ethm
Note that under the assumptions 
of (b), every parabolic subgroup 
is solvable, and every loxodromic 
subgroup  is either solvable or contains 
a non-abelian free subgroup. 

In Section \ref{sec:ind-grps} 
we collect
a necessary information about (nested) 
ind-groups of automorphisms 
of an affine variety and on subgroups 
with no finite index subgroups. 

The main result of 
Section \ref{sec:solv-conn} is 
the following theorem, see 
Theorem \ref{prop:clef} and 
Corollary \ref{cor:clef}.
\bthm\la{mthm-2} Let $X$ be an affine variety, 
$G=A*_{A\cap B} B\subset\Aut(X)$ be
an amalgam with the corresponding
Bass--Serre tree $T=T(G)$, and $H$ 
be a solvable  connected 
subgroup  of $G$. 
Assume that every $h\in H$ 
whose fixed point set $T^h$ is unbounded,
 is a torsion element. 
Then $H$ is conjugate 
to a subgroup of one of the factors 
$A$ and $B$. In particular, 
if every $g\in G$ with unbounded 
$T^g$ is a torsion element, 
then every Borel subgroup
of $G$ is conjugate to 
a Borel subgroup 
of $A$ or $B$.
\ethm
This theorem is inspired by 
some known results for
the group $\Aut(\AA^2)$, see 
\cite[Chapter 7]{Lam25}.
Our proof of Theorem \ref{mthm-2}
uses a result of Cantat, Regeta and Xie 
 on algebraicity of certain
connected commutative
groups of automorphisms, see 
 \cite[Theorem B]{CRX23} 
 or Theorem~\ref{thm:CRX} below.

In Section \ref{sec:A2} 
we provide a description 
of the Borel subgroups of 
$\Aut(\AA^2)$ 
and $\Aut_0(\AA^2)$, where 
$\Aut_0(\AA^2)$ is the 
subgroup of automorphisms 
of $\AA^2$ fixing the origin. 
Here, we mainly follow the papers 
\cite{BEE16} and \cite{FP18}. 
We also describe the centralizers of 
torsion elements according to \cite{AZ13}. 

In subsection \ref{ss:max-solv-A2}, 
we study the solvable 
(not necessarily connected)
subgroups  of $\Aut(\AA^2)$ that 
are maximal among 
the solvable subgroups.
According to Lamy's classification, 
see \cite[Chapter~7]{Lam25},
there are four types 
of subgroups  of $\Aut(\AA^2)$
based on their action 
on the Bass-Serre tree $T$
that corresponds to 
Jung--van der Kulk--Nagata's 
amalgam structure on $\Aut(\AA^2)$.
These are
elliptic, parabolic, 
and two kinds of 
loxodromic subgroups, namely
elementary loxodromic 
and general loxodromic subgroups.
All parabolic and elementary loxodromic 
 subgroups are solvable, while
 a general 
 loxodromic subgroup contains 
 a non-abelian free group
and is therefore not solvable.

Using Theorem \ref{mthm-3} 
we prove the following result, 
see Theorems \ref{thm:Borel-A2} 
and \ref{prop:max}
and Proposition~\ref{rem:case-A2}. 
We let $\TT\subset\GL(2,\kk)$ 
be the maximal torus of diagonal matrices,
$\tau\colon (x,y)\mapsto (y,x)$
be the twist, $Z\subset\GL(2,\kk)$ 
be the subgroup 
of scalar matrices, $2O\subset\SL(2,\kk)$ 
be the binary octohedral subgroup
of order 48,
and $\Trans(\AA^2)\subset\Aff(\AA^2)$ be
the subgroup of translations. 
\bthm\la{mthm-5}
Every maximal solvable subgroup 
of $\Aut(\AA^2)$ is conjugate to 
precisely one of the following maximal 
solvable subgroups of $\Aut(\AA^2)$:
\begin{itemize}
\item the triangular de Jonqui\`eres subgroup 
$\cB$;
\item
$\langle \TT, \tau, \Trans(\AA^2)\rangle
\subset\Aff(\AA^2)$; 
\item
$\langle Z, 2O,  \Trans(\AA^2)\rangle
\subset\Aff(\AA^2)$;
\item a parabolic subgroup 
$\Stab_{\Aut(\AA^2)}(b)$ for
a Freudenthal end $b\in \partial T$;
\item an elementary loxodromic subgroup 
$\Stab_{\Aut(\AA^2)}(\{b,b'\})$
for a pair of distinct Freudenthal 
ends $b,b'\in \partial T$.
\end{itemize}
\ethm
Section \ref{sec:toric} 
contains the proofs of 
Theorems \ref{main} 
and \ref{thm:main-2},
see Proposition \ref{prop:Borel-Nde},
Theorem \ref{thm:Borel}, 
Remark \ref{rem:tau}, and 
Theorem \ref{thm:no-fin-ind}. 

In Appendix A (Section \ref{sec:tree-amalg})
we extend some results
from Section \ref{sec:maximal} to amalgams 
of a tree of groups. We also present there
examples of affine automorphism  groups 
which are
bearable, that is,  amalgams of a tree 
of ind-subgroups. 

Appendix B (Section \ref{sec:App-B}) contains 
a description of all maximal solvable subgroups 
of the affine group $\Aff(\AA^2)$; this is 
used in the proof of Theorem \ref{mthm-5}.
\medskip

\section{Preliminaries on 
amalgams}\la{sec:prelim}
\subsection{An amalgam and 
its Bass--Serre tree}\la{ss:BS-tree}
In this subsection, we mainly follow  
\cite{Bas76}, \cite[Appendix D]{Lam25}, 
\cite{Ser03}, \cite{Tit70}, \cite{Tit77}, 
and \cite{Wri92}. We notably use
the terminology of \cite{Lam25}.

Consider an amalgam $G=A*_C B$ 
of two subgroups $A, B\subset G$
sharing a common subgroup 
$C=A\cap B$, see e.g. \cite{Ser03}. 
If $C$ is proper in $A$ and in $B$, 
then we say that the amalgam 
$G=A*_C B$ is proper. 
We assume by default that all amalgams 
considered in this article are proper.
Given an amalgam $G=A*_C B$, 
the associated \emph{Bass--Serre tree} 
$T=T(G)$ 
is a bicolored tree 
whose vertices 
are the left cosets $G/A$ and $G/B$, 
and whose edges are the left cosets $G/C$.
The vertices $gA$ and $gB$ for 
$g\in G$ 
are joined by the edge $gC$. 
The vertices $gA$ and $g'B$ 
are joined by an edge if and only if the 
left cosets $gA$ and $g'B$ 
intersect, and in this case 
$gA\cap g'B=g''C$.
Note that the tree $T$ is not 
necessarily locally finite.

The group $G$ acts on $T$ 
via the left multiplication, 
for example $f\in G$ 
acts via $f.gA=(fg)A$, etc.
It acts transitively on edges 
without inversions, that is,
if an edge $gC=[gA,gB]$
is fixed by $f$, then the vertices 
$gA$ and $gB$  are also fixed by $f$. 
Thus, $G$ acts transitively 
on vertices of type $gA$
and on vertices of type $gB$, and does 
not mix them.
The stabilizer of a vertex $gA$ (resp. $gB$) 
is the conjugate subgroup $A^g:=gAg^{-1}$ 
(resp. $B^g=gBg^{-1}$),
and the stabilizer of an edge $gC$ 
is the conjugate
subgroup $C^g=gCg^{-1}$. Thus, $f$ 
fixes a vertex (resp. an edge) 
of $T$ 
if and only $f$ is conjugate to 
an element of $A\cup B$ 
(resp. of $C$). 

 It is worth to note that each 
of the subgroups $A$ and $B$
coincides with its normalizer in $G$. 
In other words, $A$ resp. $B$ 
acting on $T$ fixes 
a unique vertex of $T$
denoted by $1\cdot A$ resp. $1\cdot B$ 
(and $C$ fixes the
edge $1\cdot C=[1\cdot A-1\cdot B]$
of $T$).
The stabilizers of vertices 
are pairwise distinct.
\bdefi[\emph{induced amalgam structure}]
\la{induced-amalg}
Let $G=A*_C B$ be an amalgam, where 
${C=A\cap B}$ is a proper subgroup of $A$ 
and of $B$. Suppose that subgroups 
$A_0\subset A$ and $B_0\subset B$
satisfy $A_0\cap C=B_0\cap C=:C_0$, 
where $C_0$ is proper in $A_0$ and in $B_0$.
Then $G_0:=\langle A_0, B_0\rangle\subset G$
inherits an amalgam structure 
$G_0=A_0*_{C_0} B_0$, 
see \cite[Corollary~8.11]{Neu48}.
\edefi
\brem\la{rem:tree-embed}
Due to the uniqueness of the 
reduced word presentation,
we have 
\[A_0=A\cap G_0,\quad B_0=B\cap G_0,
\quad\text{and}\quad C_0=C\cap G_0.\]
Therefore, the correspondence
\[g_0 A_0\mapsto g_0 A, \quad
g_0 B_0\mapsto g_0 B,\quad
g_0 C_0\mapsto g_0 C\quad\forall g_0\in G_0\]
defines an embedding of trees 
$T(G_0)\hookrightarrow T(G)$.
\erem
\bdefis\la{el-lox}
The 
\emph{translation length} of $f\in\Aut(T)$ is
 \[{\rm tl}(f)=\min\{{\dist}_T(v, f(v))\,|\, v\in 
 \vert(T)\},\]
 where $\dist_T$ stands for the tree distance on $T$.
An element $f\in G$ is said to be \emph{elliptic} 
if  ${\rm tl}(f)=0$, i.e. if $f$ fixes a vertex of $T$,
and \emph{loxodromic} otherwise 
\footnote{Alternatively, 
a loxodromic element is said to be 
\emph{hyperbolic} or \emph{axial}, 
see e.g. \cite{Ser03}, 
\cite{FL10}, and \cite{Syk02}.}.
\edefis
\blem[{\rm see \cite[Proposition~25]{Ser03} 
or \cite[Lemma D.8]{Lam25}}] 
\la{lem:loxo-elem}
An element $g\in G$ is loxodromic 
if and only if 
$g$ is conjugate in $G$ to a
cyclically reduced element of even length $2n$,
and then ${\rm  tl}(g)=2n$. 
\elem
\brems\la{rem:loxod} $\,$

1. Some of the results we cite
hold more generally for 
groups acting on trees. Indeed, 
let $G$ be a group 
acting on a tree $T$ without inversions 
and let $e=[u,v]$ be an edge of $T$ 
connecting the vertices $u$ and $v$ of $T$. 
Let  $G_u, G_v$ and $G_e$ be the 
stabilizers of $u$, $v$ and $e$ in $G$.
Then $G_u*_{G_e} G_v$ 
is naturally isomorphic to
the subgroup $\langle G_u, G_v\rangle$
of $G$, see 
\cite[Corollaire~4.3]{Hau81}. 
If, moreover, $G$ 
acts faithfully and transitively
on the set of edges of $T$, 
then 
$G\simeq G_u*_{G_e} G_v$,
 see \cite[Corollaire~4.4]{Hau81}.

2. 
According to 
\cite[Proposition~3.2]{Tit70} 
every automorphism $f$ of $T$ 
acting without inversions 
either 
fixes a vertex (and we then say 
that $f$ is elliptic), 
or leaves invariant a unique 
infinite geodesic $\Ax(f)$,
 that is, a path without backtracking
infinite in both directions, on 
which it acts as a shift 
(and we then say that $f$ is loxodromic). 
The fixed subgraph $T^f$ 
of an elliptic $f$ is a subtree of $T$, 
see \cite[Section~6.1]{Ser03}. 
A loxodromic automorphism $g$ of $T$
acts via a shift by ${\rm  tl}(g)$ on 
the unique invariant geodesic $\Ax(g)$,
and
\[\Ax(g) := \{x \in \vert(T)\,|\, 
{\rm dist}_T(x,g(x)) ={\rm  tl}(g)\},\]
see  \cite[Proposition~24]{Ser03},
\cite[Proposition~3.2]{Tit70}, 
or \cite[Proposition~3.3]{Hau81}. Thus, 
 $\langle g\rangle\simeq\ZZ$.
 \erems
\bprop \la{prop:elements}
Let $f, f_1,f_2\in G=A*_C B$ be elliptic, 
and $g, g_1, g_2\in G$ be loxodromic. 
Then the following hold.
\bnum
\item[$(a)$] The product $f_1f_2$ 
is loxodromic 
if and only if 
$T^{f_1}$ and $T^{f_2}$ are disjoint. 

\item[$(b)$]  
If $f$ and $g$ commute, then 
$T^f\supset {\rm Ax}(g)$, 
and therefore
$T^f$ is unbounded, that is, has infinite 
diameter in the tree metric.
\item[$(c)$]
Let $e$  be an edge of ${\rm Ax}(g)$, $G_e$ 
the stabilizer of $e$ in $G$, and 
$F=G_e\cap \Cent_G(g)$.
Then $F$ is a normal subgroup 
of $\Cent_G(g)$ and 
$\Cent_G(g)/F\simeq\ZZ$. 
\item[$(d)$] An element
$h\in G$ is loxodromic 
if and only 
if, for some $v\in \vert(T)$, the orbit $
\langle h\rangle (v)$ is unbounded, 
if and only if every  orbit of $
\langle h\rangle$ acting on $\vert(T)$ 
is unbounded.
\item[$(e)$] The axis $\Ax(g)$
of a loxodromic $g\in G$ 
 is contained in every subtree 
$T'\subset T$ 
stable under $g$ and $g^{-1}$, 
and therefore $g|_{T'}$ is loxodromic.
\item[$(f)$] If $\Ax(g_1)\cap\Ax(g_2)$ 
is finite, then the subgroup 
$\langle g_1, g_2\rangle$
 contains a free group of rank $2$.
\item[$(g)$] A subgroup $H\subset G$ 
acting freely on $T$ is a free group. 
\item[$(h)$] Let $H\subset G$ be 
a subgroup and let 
$N\subset H$ be the subgroup 
generated by the
stabilizers $H_v$ of vertices 
$v\in \vert(T)$.
Then $N$ is a normal subgroup 
of $H$ and $H/N$ 
is a free group.
\enum
\eprop
\bproof 
See e.g. \cite[Lemme~2]{PV91} 
and the references therein for (a) and
\cite[Lemme~3.1]{Lam01} for (b).
See also
\cite[Lemma~2.1]{Syk02} for (c), 
\cite[Corollaire~3.3 and the proof 
of Corollaire~3.5]{Tit70} for (d) and (e),
 \cite[Lemme~2.3]{Hau81}, 
 \cite[Lemma~2.6]{CM87} or
\cite[Lemme~6 and Proposition~3]{PV91}
 for (f),
\cite[Proposition~18]{Ser03} for (g), and 
\cite[Corollary~15]{RS} for (h).
\eproof
\bdefis\la{def:end} $\,$

\begin{itemize}
\item Cutting an edge $[u,v]$ of 
an infinite geodesic $\gamma$ 
produces two 
half-geodesics, also  called 
\emph{geodesic rays}. 

\item A \emph{boundary point}, 
or an \emph{end} 
$b\in\partial T$ is
the equivalence class of geodesic rays
 in $T$, 
where two geodesic rays are 
equivalent if their 
intersection contains a geodesic ray. 
\end{itemize}
\edefis

\brems\la{rem:ends} $\,$

1. We let $\partial T$ 
be the set of ends of $T$.
For each end $b\in \partial T$ 
and each vertex $v$ of $T$
there is a unique geodesic ray 
$r=:[v,b[\in b$
with extremity (or tip) $v\in\vert(T)$. 
Every infinite geodesic $\gamma$ of $T$
belongs to 
exactly two distinct ends 
$b,b'\in \partial T$; we  then write
$\gamma=]b,b'[$. For every pair 
$b,b'\in \partial T$ there exists a 
unique infinite geodesic $\gamma$
of $T$ such that $\gamma=]b,b'[$, 
see \cite[Section~1.5]{Tit77}. 

2. The action of $\Aut(T)$ 
on $T$ naturally extends 
to the boundary $\partial T$.
One says that $f\in\Aut(T)$ 
\emph{fixes an end $b\in\partial T$}
if for some (and then for any) geodesic ray
$r\in b$ the intersection $r\cap f(r)$ 
contains a geodesic ray $r'\in b$. 
\erems

\blem[{\rm see e.g. \cite[p.~157]{PV91}}] 
\la{lem:end}
Let $f\in\Aut(T)$ be elliptic and $v\in T^f$. 
If $f$ fixes an end $b\in\partial T$, 
then it pointwise fixes  the geodesic ray 
$r=[v,b[\in b$ with a tip $v$.
A loxodromic $g\in\Aut(T)$ fixes $b$ 
if and only if $\Ax(g)$
contains a geodesic ray $r\in b$. 
\elem
\subsection{Subgroups of an amalgam}
\la{ss:sbgrps}
\bdefis\la{def:sbgrps} $\,$

1.
A subgroup $H\subset G=A*_C B$ 
that fixes 
a vertex of $T$ is called 
\emph{elliptic}. Thus $H$ is elliptic 
if and only if 
it is conjugate to a subgroup of 
$A$ or $B$.
An elliptic subgroup $H$ is made up
of elliptic elements.

2. 
However, a subgroup $H$ made 
up of elliptic 
elements is not necessarily elliptic.
If such a subgroup $H$  is not elliptic, 
it is called
\emph{parabolic}. 

3. 
A subgroup $H$ 
is said to be \emph{loxodromic} if 
it contains a 
loxodromic element. 
A loxodromic subgroup $H$ is said to be
\begin{itemize}
\item \emph{elementary loxodromic} 
if $H$ contains a loxodromic element $g$
such that all elements
in $H$ preserve the axis ${\rm Ax}(g)$;
\item \emph{focal loxodromic} 
if $H$ fixes a
unique boundary point  $b\in \partial T$ 
\footnote{Cf. Proposition \ref{thm:Tits}(d).};
\item \emph{general loxodromic} otherwise.
\end{itemize}
\edefis
Another classification of groups acting on 
trees is presented in 
\cite[Sections~2 and~3]{Tit77};
 cf. also \cite{KS70}. 
For different types of
 subgroups of an amalgam,
 we have the following results. 
\bprop\la{thm:Tits} 
Let $G=A*_{A\cap B} B$
be a proper amalgam 
and $T=T(G)$ be its Bass-Serre tree. 
\begin{itemize}
\item[$(a)$]  
A solvable torsion subgroup $H\subset G$ 
 is either elliptic or parabolic. 
\item[$(b)$] A solvable loxodromic 
subgroup $H\subset G$
leaves invariant an  end  
or a pair of  ends
of $T$. 
\item[$(c)$] A parabolic subgroup 
$H$ of $G$
fixes a unique end
$b\in\partial T$ and is not 
finitely generated. 
\item[$(d)$] An elementary
 loxodromic subgroup $H$
leaves invariant a unique pair 
$\{b,b'\}$ of distinct ends 
of $T$.
\item[$(e)$] A general loxodromic 
 subgroup $H\subset G$
contains a free group $F$ 
of rank $2$ generated 
by two loxodromic elements 
with disjoint axes. 
The group $F$ acts freely on $T$ 
and fixes no
point of $\partial T$.
\end{itemize}
\eprop
\bproof 
Any loxodromic element of $G$ 
is of infinite order, 
see Remark \ref{rem:loxod}.2. 
This implies~(a). 
Statement (b) follows e.g. from 
\cite[Corollary~2]{Tit77}.
See \cite[Corollary~3 of 
Proposition~26]{Ser03}
and \cite[Proposition~3.4]{Tit70}; 
cf. also \cite[Corollary~3]{Tit77} 
and \cite[Lemma~D.11]{Lam25} for (c). 

Statement (d) is true, because
otherwise $H$ would have to
leave invariant 
at least three geodesic axis 
$]b_i, b_j[$ for $i,j=1,2,3$, $i\neq j$, 
which is impossible by
Proposition~\ref{prop:elements}(e). 
See \cite[Proposition~3.4 and~3.5]{Hau81}
and \cite[Proposition~4.3]{Lam01} for (e).
\eproof
\brems\la{rem:classification}
Consider a not necessary 
solvable subgroup $H$ of 
a proper amalgam ${G=A*_{C} B}$
that fixes no vertex of $T=T(G)$, 
 that is, is not elliptic. 
According to the classification of 
Definition \ref{def:sbgrps}, $H$ has the 
following properties. 
\bnum
 \item[$1.$] 
A parabolic subgroup 
 $H=\cup_{i} H_i$
 is filtered by a non-stationary 
 ascending sequence 
 of subgroups. Every $H_i$ is conjugate 
 to a subgroup of $C$ and fixes 
 a geodesic ray 
 $r_i\in b$, where $r_i\supset r_{i+1}$ 
 and $\cap_i r_i=\emptyset$, 
cf. \cite{Mol67}, \cite{Nag80}, and 
\cite[Theorem~0.3(2)]{Wri79}. 
 
\item[$2.$] 
An abelian loxodromic subgroup $H$
is a product $F\times \langle g\rangle$, 
where $F$ is conjugate to a subgroup of
$C$ and $g$ is loxodromic, see 
\cite[Theorem~0.3(3)]{Wri79}. 
See also \cite[Section~2]{Wri79} 
and \cite[Examples~7.65--7.66]{Lam25}
for examples
of abelian parabolic resp. loxodromic subgroups
of an amalgam. 
The parabolic subgroups in these examples 
are isomorphic to 
subgroups of the group $\cR\subset\Gm$ 
of all roots of unity.

\item[$3.$] Let $H\subset G$ 
be a focal loxodromic subgroup
fixing an end $b\in\partial T$.
From Remark ~\ref{rem:loxod}, 
Definitions \ref{def:end} and \ref{def:sbgrps}, 
and Lemma \ref{lem:end} 
we extract the following.
The fixed subtrees $T^{f_i}$ of each pair 
of elliptic elements $(f_1,f_2)\in H^2$ share 
a common geodesic ray representing $b$.
The same holds for the axes $\Ax(g_i)$ 
of each pair 
of loxodromic elements $(g_1,g_2)$
 and for 
$\Ax(g)$ and $T^f$, where
$(g,f)\in H^2$ with 
$g$ loxodromic and $f$ elliptic.
Moreover, there exists a pair of 
loxodromic elements
$(g_1,g_2)\in H^2$  with 
$\Ax(g_1)\neq\Ax(g_2)$. 
See \cite[D10]{Lam25}
for examples of focal loxodromic subgroups. 

 \item[$4.$] 
 A subgroup $H\subset G$ that 
does not fix any point of $T\cup\partial T$ 
is loxodromic, see  \cite[Proposition~1]{PV91}. 
If $H$ with $T^H=\emptyset$
leaves no point and no pair  
of points on $\partial T$ invariant, 
then it is general loxodromic, 
see  \cite[Proposition~2]{PV91}.
\item[$5.$] 
We note that some of the above results
are also valid for groups acting 
on general graphs, 
see e.g. \cite{Jun94} 
and the references therein.
\enum
\erems
\bdefi
A subset $S$ of $G=A*_{C} B$ 
is said to be \emph{of bounded length} 
if there is an upper bound 
on the lengths 
 of reduced decompositions 
 of elements of $S$.
\edefi
The following theorem is due to Serre  
\cite[Theorem~8]{Ser03}. 
 \bthm
 \la{thm:Serre}
Every subgroup 
$H\subset G=A*_{C} B$ 
 of bounded length
is elliptic, that is, is contained in a conjugate 
of $A$ or $B$.
\ethm
\section{Maximal solvable subgroups 
of amalgams} \la{sec:maximal}
In this
section, we present some
classes of
maximal solvable subgroups of 
amalgamated free products.
\bprop\label{thm:maximal}
Consider an amalgam $G=A*_C B$, 
where $C=A\cap B$, $|A:C|\ge 3$ 
and $|B:C|\ge 2$.  
If $A$ is solvable, 
then $A$ is maximal among the
solvable subgroups of $G$. 
\eprop
\bproof
Suppose that $A$ is solvable and
 there exists a solvable 
subgroup $H$ of $G$ properly containing
$A$. We can then
find an element $h\in H\setminus A$ 
expressed as a strictly alternating word 
 in $G$ of the form 
\[h=b_1a_1\cdots b_{k-1}a_{k-1}b_{k},\quad\text{
where}\quad k\ge 1, \quad a_i\in A\setminus B
\quad\text{and}\quad b_i\in B\setminus A
 \,\,\, \forall i.\] 
Thus, $h$ has length $l_G(h)=2k-1\ge 1$.
Since $|A:C|\ge 3$ we can
choose the elements $\alpha_0, \alpha_1\in 
A\setminus B$ 
from distinct right cosets modulo $C$,  that is,
$C\alpha_0\neq C\alpha_1$. 
For $ j=0,1$ consider the commutators
$f_j=[\alpha_j,h]\in H$. 
We claim that $F:=\langle f_0, f_1\rangle\subset H$ 
is a free group of rank two. Since $H$ is solvable, 
this gives a contradiction and proves our assertion.

We have 
\[f_j=\alpha_j^{-1}h^{-1}\alpha_jh
=\alpha_j^{-1}b_{k}^{-1}a_{k-1}^{-1}\cdots 
a_1^{-1}b_1^{-1}\alpha_jb_1a_1\cdots 
a_{k-1}b_{k}.\] 
Therefore, $f_j$ is cyclically reduced 
with length $l_G(f_j)=4k\ge 4$.

We assert that any 
nonempty reduced word $w(f_0,f_1)$
considered as a strictly alternating word in $G$ 
has positive length. In fact, 
the only possible cancellations in $w(f_0,f_1)$ 
occur in the syllables 
$f_0f_1^{-1}$ and $f_1f_0^{-1}$.
For example, the length in $G$ of the syllable
\[f_0f_1^{-1}=\alpha_0^{-1}h^{-1}
\alpha_0\alpha_1^{-1}h\alpha_1\]
equals $4k+1\ge 5$. Indeed, 
since $C\alpha_0\neq C\alpha_1$ 
we have 
$\alpha_0\alpha_1^{-1}\in A\setminus C
=A\setminus B$. 
Thus, after all possible cancellations 
in a reduced word 
$w(f_0,f_1)\in G$ 
of positive length in the free group 
with generators 
$f_0$ and $f_1$, its length in $G$ 
is also positive. 
Hence 
$F=\langle f_0, f_1\rangle\subset H$ 
is a free group of rank two. 
\eproof
\brems\la{rem:bearable} $\,$

1. 
The conclusion of Proposition 
\ref{thm:maximal} 
is no longer valid in the case where 
$|A:C|=2$ and $|B:C|=2$.  
As an example, 
consider the infinite 
dihedral group
\[D_\infty=\ZZ/2\ZZ*\ZZ/2\ZZ\simeq 
\ZZ\rtimes (\ZZ/2\ZZ).\] 
It is solvable, therefore
the factors of the free product 
are not maximal solvable subgroups. 
The Bass-Serre tree $T=T(G)$
is a chain; it is 
a unique geodesic of $T$.
The factor $\ZZ$ of the 
semi-direct product
acts 
on $T$ by translations, and
the generator of $\ZZ/2\ZZ$
acts by symmetry. Thus, $D_\infty$
is a non-proper elementary loxodromic 
subgroup of itself. 

2.
Suppose that, under the assumptions of 
Proposition \ref{thm:maximal}, $B$ is 
solvable, $|B:C|=2$ and $|A:C|\ge 3$. 
We do not know any example where
such a subgroup $B$ is not maximal 
among the solvable subgroups of $G$. 
For example, in \eqref{eq-amalg-2} 
we consider an amalgam 
$\cN_{d,e}=\cN_{d,e}^+*_{N_{d,e}^+} N_{d,e}$,
where $\cN_{d,e}^+\subset \cN_{d,e}$
are subgroups of 
$\Aut_0(\AA^2)=\{f\in\Aut(\AA^2)\,|\,f(0)=0\}$,
see Section \ref{sec:toric} for notation.
The positive integers $d$ and $e$ satisfy
$\gcd(d,e)=1$ and $e^2\equiv 1\mod d$. 
For $e\not\equiv 1\mod d$, the subgroup 
$N^+_{d,e}=\TT$
is the  maximal torus of $\GL(2,\kk)$ consisting 
of the diagonal matrices. We have
$N_{d,e}\simeq\TT\rtimes \ZZ/2\ZZ$. 
The subgroups $A=\cN^+_{d,e}$, $B=N_{d,e}$, 
and $C=A\cap B=\TT$
satisfy the  above conditions. 
By Proposition \ref{prop:max}(a), $B$ 
with $|B:C|=2$ is a maximal 
solvable subgroup of $\Aut_0(\AA^2)$.
Therefore, it also is a maximal 
solvable subgroup of 
$A*_{C} B=\cN_{d,e}\subset \Aut_0(\AA^2)$.
\erems
The maximal solvable
subgroups as in 
Proposition \ref{thm:maximal} are elliptic.
In Proposition \ref{prop:exas-max-loxo} 
and Corollary \ref{cor:lox-solv2}
we describe, for certain amalgams, 
parabolic and  loxodromic 
subgroups that are maximal among 
the solvable subgroups. 
In the case of the 
Jung-van der Kulk-Nagata 
amalgam structure
on $\Aut(\AA^2)$, 
some examples of this type 
come from
\cite[Proposition~4.10]{Lam01}, 
cf. also
\cite[Section~5.1, Remark~2]{BEE16}.

In the rest of this section we 
use the following notation.
\bnota
Consider an amalgam $G=A*_C B$
acting on 
its Bass-Serre tree $T=T(G)$.
Given an end $b\in\partial T$,
we denote by $S_b=\Stab_G(b)$
the subgroup of all $g\in G$ fixing $b$
and by $P_b\subset S_b$ the subgroup of all 
elliptic elements of $S_b$. Given a loxodromic $g\in G$,
we denote by $H(g)$ the subgroup of all elements of $G$ 
 leaving $\Ax(g)$ invariant. 
Thus, $H$ is a maximal elementary loxodromic 
subgroup, see Definition \ref{def:sbgrps}. 
\enota
\blem\la{lem:loxo-solv} 
Suppose that the subgroup
$C=A\cap B$ is solvable, with derived length $l$. 
Then, 
the subgroups $S_b$ and $P_b$
are solvable, with derived lengths $\le l+1$
and $\le l$, respectively. Moreover, 
 all parabolic, elementary loxodromic, and 
focal loxodromic subgroups of $G$
are solvable with derived length $\le l+1$. 
\elem
\bproof Choose a ray $r=[v,b[\in b$, and let 
$v_0=v,v_1,\ldots$ be the sequence of 
vertices of $r$ enumerated such that 
$\dist_T(v,v_i)=i$. Let $S_i$ 
be the subgroup of all $s\in S_b$ which fix 
the ray $r_i:=[v_i, b[\in b$.
It is easy to see that  for any pair
$s,s'\in S_b$, there exists $i\in\NN$ 
such that $[s,s']\in S_i$. Hence, 
$S_b^{(1)}=[S_b,S_b]\subset \bigcup_i S_i$.

The subgroup $S_i$
fixes the edge $[v_i, v_{i+1}]$ of $r$.
 It is therefore
conjugate to a subgroup of~$C$. 
It follows that $S_i^{(l)}=1$ for all $i$.
So also 
$S_b^{(l+1)}\subset \bigcup_i S_i^{(l)}=1$.
Thus, $S_b$ is solvable with derived length 
$\le l+1$. Since $P_b\subset \bigcup_i S_i$,
it follows that $P_b$ is solvable 
with derived length $\le l$. 

Note that a parabolic subgroup $P$ fixing 
an end $b\in\partial T$
is contained in $P_b$ 
and a focal loxodromic subgroup $S$ fixing 
an end $b\in\partial T$
is contained in $S_b$,
see
 Definitions \ref{def:sbgrps}
 and Proposition \ref{thm:Tits}(c). 
If $H$ is 
an elementary loxodromic subgroup 
with axis $]b,b'[$, then the subgroup 
$H_0\subset H$ consisting of all 
$h\in H$ which fix
$b$ and $b'$ is contained in 
$S_b\cap S_{b'}$ and has index 
$\le 2$ in $H$. 
Now the last assertion follows. 
\eproof
\bprop\la{prop:exas-max-loxo}
Consider a proper amalgam $G=A*_C B$
acting on 
the Bass-Serre tree $T=T(G)$.
Suppose that the subgroup
$C=A\cap B$ is solvable. 
Let $b\in\partial T$ be such that 
$S_b$ is not elliptic. 
Then one of the following statements is true.
\begin{itemize}
\item[$(i)$] $S_b=P_b$ is parabolic 
and maximal 
among the solvable subgroups of $G$;
\item[$(ii)$]  $S_b$ is focal loxodromic 
and maximal 
among the solvable subgroups of $G$;
\item[$(iii)$] $S_b$ is 
elementary loxodromic
and is contained in 
a maximal solvable subgroup $H$,
where either $H=H(g)$ is 
a maximal elementary loxodromic
subgroup containing 
a loxodromic $g\in S_b$, 
or $H=S_{b'}$ is 
a focal loxodromic subgroup, where
$\Ax(g)=]b,b'[$. 
\end{itemize}
\eprop
\bproof 
Suppose first that $S_b=P_b$ is parabolic.
Then $b$ is the unique end fixed by $S_b$,
see Proposition \ref{thm:Tits}(c).  
Suppose that $R\supset S_b$ is 
a solvable subgroup of 
$G$ strictly containing $S_b$.
Then $R$ cannot fix $b$ and is neither elliptic, 
nor parabolic, nor general loxodromic. 

If $R$ is elementary loxodromic, then $R$ 
leaves invariant 
the pair of ends $(b,b')$ of the axis 
$\Ax(R)=]b,b'[$ of $R$. 
Since $S_b\subset R$ fixes $b$, 
it fixes also the second end 
$b'$, which contradicts the uniqueness 
of the end fixed by $S_b$.

If $R$ is focal loxodromic, 
then it fixes a unique end 
$\hat b\neq b$, and therefore 
$S_b\subset R$ 
fixes $b$ and~$\hat b$, 
which is impossible. 
This proves that  $S_b$ 
is maximal solvable, 
as indicated in (i). 

Now, suppose that $S_b$ is focal loxodromic.
Assume on the contrary that $R\subset G$ 
is a solvable subgroup 
strictly containing $S_b$. Then $R$ 
cannot fix $b$ and is neither elliptic, 
nor parabolic, nor elementary loxodromic, 
nor general loxodromic. 
Being focal loxodromic, 
$R$ fixes a unique end 
$\hat b\neq b$ of $T$, 
and therefore  $S_b\subset R$ 
fixes two 
distinct ends $b$ and $\hat b$, 
which is impossible. In this case,
 $S_b$ is maximal solvable, 
 as stated in (ii). 
 
 Finally, assume that $S_b$ is 
 elementary loxodromic.
If $g\in S$ is loxodromic, 
then $S_b\subset H(g)$. 
Let again $\Ax(g)=]b,b'[$, 
where $b'\in\partial T$, 
$b'\neq b$.
 If $S_{b'}\not\subset H(g)$, then $S_{b'}$ 
is neither elementary  loxodromic 
nor general loxodromic,
see Lemma \ref{lem:loxo-solv}.
 Therefore, it is focal loxodromic 
 and maximal solvable
 thanks to (ii). Since $S_b\subset H(g)$ 
 fixes $b$, it also fixes $b'$, 
 and so $S_b\subset S_{b'}=:H$.
 
 If  $S_{b'}\subset H(g)$, 
 then $S_{b'}=S_b$ 
 is a subgroup of index at most 2 in $H(g)$. 
We assert that, in this case, $H(g)$ 
is maximal solvable. Indeed,
suppose $H(g)\subset R$, 
where $R$ is a solvable
loxodromic subgroup of $G$. 

If $R$ is 
elementary loxodromic, then $R=H(g)$ 
due to the maximality of $H(g)$ among 
the elementary loxodromic 
subgroups of $G$
that contain $g$. 

Assume that $R$ 
strictly contains $H(g)$. 
Being solvable, $R$ is 
not general loxodromic. 
It is therefore focal loxodromic
and fixes a unique end 
$\hat b\in\partial T$. Then $\hat b$  is 
also fixed by  $H(g)$. 
It follows that the geodesic $]b,\hat b[$ 
is invariant under $H(g)$, 
and therefore coincides with $]b,b'[=\Ax(g)$.
Thus, $R$ fixes $\hat b\in\{b,b'\}$, 
hence coincides with the subgroup 
$S_b=S_{b'}\subset H(g)$. This contradiction 
proves that $H(g)$ is maximal solvable. 
\eproof
The following corollary is immediate. 
It gives the conditions necessary  for a subgroup 
of an amalgam to be maximal solvable. 
\bcor\la{cor:lox-solv1} Let the assumptions 
of Proposition \ref{prop:exas-max-loxo} 
be verified. 
Let  $H$ be a maximal solvable subgroup 
of $G$.
Then precisely one of the 
following statements hold.
\begin{itemize}
\item[$(i)$] $H$ is elliptic and conjugate 
to a maximal solvable subgroup of 
$A$ or $B$;
\item[$(ii)$]  $H$ is parabolic and 
$H=P_b=S_b$ 
for some $b\in\partial T$;
\item[$(iii)$]  $H$ is elementary 
loxodromic  and $H=H(g)$
for a loxodromic $g\in G$;
\item[$(iv)$]  $H$ is focal loxodromic 
and $H=S_b$ 
for some $b\in\partial T$.
\end{itemize}
\ecor
We have the 
following criteria 
for a subgroup of an amalgam 
to be maximal solvable. 
\bcor\la{cor:lox-solv2} 
Under the assumptions 
of Proposition \ref{prop:exas-max-loxo}, 
let us further assume that $G$ has 
no focal loxodromic subgroup. 
Let $H\subset G$ be a non-elliptic subgroup. 
Then $H$
is maximal solvable if and only if 
one of the following conditions holds:
\begin{itemize}
\item[$(i)$]  $H$ is parabolic and 
$H=P_b=S_b$; 
\item[$(ii)$]  $H$ is elementary 
loxodromic  and $H=H(g)$
for a loxodromic $g\in G$.
\end{itemize}
\ecor
\bproof By Lemma \ref{lem:loxo-solv}, any $H$
as in (i) and (ii) is solvable. 
By Corollary 
\ref{cor:lox-solv1}(ii) resp. (iii),
the conditions (i) resp. (ii) are necessary. 
To prove that they are sufficient, suppose
first that (i) holds. 
It suffices to show that
$H$ is not contained in any subgroup $H(g)$,
where $g$ is loxodromic. Indeed, if 
$H=S_b\subset H(g)$, then $S_b$ 
fixes the end
$b'\neq b$ of $\Ax(g)=]b,b'[$. 
This gives a contradiction, since $b\in\partial T$
is a unique end fixed by the parabolic
subgroup $H=P_b$.

Now suppose that (ii) holds. 
Then $H$ is maximal among 
elementary loxodromic subgroups of $G$,
 and is not contained 
in a  focal  loxodromic or 
general loxodromic solvable 
subgroup, 
because there is no such subgroup of $G$. 
The proof is complete. 
\eproof
\brem
The maximal elementary loxodromic 
subgroup $H(g)$ is 
not always proper
in~$G$. For example, $H(g)=G=D_\infty$ for
the amalgam of 
Remark \ref{rem:bearable}.1.
In this example, $\langle g\rangle\simeq \ZZ$ 
is a proper subgroup
of $H(g)$ of index 2 which
coincides with the centralizer of $g$;
cf. \cite[Proposition~4.8]{Lam01}.

In fact, \emph{$H(g)=G$ if and only if $C$ 
is a subgroup 
of index $2$ in $A$ and in $B$. }

Indeed, 
assume that $H(g)=G$. 
Then $H(g)$ contains $A$ and $B$, 
which therefore 
preserve $\Ax(g)$. Up to conjugation, 
we can assume that $\Ax(g)$ 
contains the edge 
$1\cdot C=[1\cdot A-1\cdot B]$.
Then $A$ and $B$ act on $\Ax(g)$ by
fixing the vertices 
$1\cdot A$ and $1\cdot B$ 
of $\Ax(g)$,  respectively, 
and $C$ fixes them  both. 
Thus, $C$ acts
on $\Ax(g)$ identically, and the  
elements of $(A\cup B)\setminus C$ act 
of $\Ax(g)$ by central symmetries. 
It follows that $C$ is a common subgroup 
of index $2$ in $A$ and in $B$. 

As for the converse implication,  
suppose that $|A:C|=|B:C|=2$. 
Then  $T=T(G)$ is an infinite 
chain, which is the only 
geodesic of $T$. 
Therefore, $G=H(g)$. 
\erem
\section{Ind-groups and groups 
without finite-index 
proper subgroups}\la{sec:ind-grps}
\subsection{Ind-groups of automorphisms}
\la{ss:ind-grps}
Recall the following definitions; 
see e.g. \cite{FK18} or \cite{KPZ17}.
\bdefis $\,$
\begin{itemize}
\item
An \emph{ind-variety} 
$\fG=\varinjlim_{k\in\NN} \fG_k$ is 
an inductive limit of an increasing sequence 
of algebraic varieties 
$\fG_k$ with closed embeddings 
$\fG_k\subset\fG_{k+1}$. 
The \emph{Zariski ind-topology} on $\fG$ 
is defined as follows: a subset 
$\omega\subset\fG$ is closed (resp. open) 
if and only if 
$\omega\cap \fG_k$ is closed  
(resp. open) in $\fG_k$ 
for every $k\in\NN$. 

\item
An \emph{ind-group} 
is an ind-variety equipped 
with a group structure such that 
the group operation 
and the inversion are continuous 
in the ind-topology. 
An ind-group $G~=~\varinjlim_{k\in\NN} A_k$ 
is said to be \emph{affine} if the $A_k$ 
are affine algebraic varieties.

\item
A subset $M$ of an ind-variety $\fG$ 
is said to be
\emph{path connected} 
if for any two points $P,Q\in M$ there exists
$k\in\NN$ and an irreducible algebraic curve 
$C\subset M\cap \fG_k$ such that $P,Q\in C$.

\item
The \emph{connected component 
of unity $G^0$} 
is the maximal connected ind-subgroup of  
$G$; it is open and closed in $G$. 

\item
A subgroup $H$ of an ind-group $G$ 
is said to be
\emph{nested} if $H$ is filtered by 
a countable ascending chain of 
algebraic subgroups.
\end{itemize}
\edefis
\brems\la{rem:connectedness} $\,$
Let $X$  be an affine algebraic 
variety over $\kk$.

1. 
The automorphism 
group $\Aut(X)$ has a natural  
structure of an affine ind-group, 
see for example \cite{FK18} or \cite{KPZ17}. 

2. Any algebraic subgroup of an ind-group 
 is closed, see \cite{FK18}. 
 Moreover, a connected nested subgroup 
 $G \subset \Aut(X)$ is closed, 
 see \cite[Theorem~6.6]{Per23}.
 
 3. A closed subgroup of 
 a nested ind-group is nested, 
 see \cite[Remark~2.8]{KPZ17}.
 
 4. A path connected subset
 of an ind-variety $\fG$ is connected, 
see e.g. \cite[Lemma~5]{BEE16}. 
For an ind-group $G$, path connectedness 
and  connectedness are equivalent properties, 
see \cite[Remark~ 2.2.3]{FK18}.
\erems
\bthm[{\rm \cite[Corollary~2.12]{KPZ17}}]
Any nested  connected subgroup 
$G \subset \Aut(X)$ 
is of the form 
$G = R_{\mathrm u}(G)\rtimes H$, 
where $H$ is a connected 
reductive algebraic subgroup 
and $R_{\mathrm u}(G)$ is 
a nested unipotent subgroup
called the \emph{unipotent radical} of $G$. 
\ethm
\subsection{Divisible groups and groups without 
finite-index proper subgroups}
\bdefi
Recall that a group $G$ is said to 
be \emph{divisible} if 
any $g\in G$ is divisible, that is, for any $n\in\NN$ 
there exists
$f\in G$ such that $g=f^n$. 
We say that $g\in G$  is \emph{pseudo-divisible}
if there exists 
$k\in\NN$ such that $g^k$ is divisible. 
 If  every $g\in G$ is pseudo-divisible, 
 then we say that $G$ is 
\emph{pseudo-divisible}.
\edefi
Clearly, the quotient group of a (pseudo-)divisible 
group is also (pseudo-)divisible. 
\blem\la{lem:divisible} $\,$

\bnum
\item[$(a)$]
 Let $G$ be a group. 
Then every subgroup of finite index 
$H\subset G$ contains a 
normal subgroup $N\subset G$ 
of finite index.
\item[$(b)$] 
 A divisible group contains 
 no proper subgroup 
 of finite index.
 \enum
 \elem
 \bproof Statement (a) is well known, 
 see for example \cite{Wiki}. 
 As for (b), suppose $G$ is divisible 
 and contains a 
 proper normal subgroup $N\subset G$ 
 of finite index. 
 Let $g\in G\setminus N$ and 
 $p\colon G\to G/N$ be 
the natural surjection. Since $p(g)$ 
is divisible in $G/N$,
we have  $p(g)=1$, and therefore 
$g\in N$, 
which is a contradiction. 
 \eproof
 \brems $\,$
 
 1. 
 A compact topological group is divisible 
 if and only if it is connected, 
 see  \cite[Theorem~8]{Pla66}. 
 This is generally not true for 
 locally compact groups. 
For example, the connected algebraic group 
 $\SL(2,\kk)$ is not divisible. Indeed, 
 $g~=~\left(\begin{matrix}-1 & 1 \\ 
0 & -1 \end{matrix}\right)$ 
does not admit a square root in 
$\SL(2,\kk)$, 
see \cite[P.~1104]{Pla66}.

2. 
A torsion element is pseudo-divisible. 
Therefore, a torsion group
(in particular, a finite group) 
is pseudo-divisible. 

3. No nonunit element $g\in\GL(n,\ZZ)$
 is divisible,
see \cite[P.~223, Lemma]{KM79}. 
Moreover, $g$ is pseudo-divisible in 
$\GL(n,\ZZ)$
if and only if $g$ is a torsion element,
see \cite[Lemma~3.7]{LRU23}.

4.
Any algebraic group over $\kk$ 
is pseudo-divisible, see 
\cite[Lemma~3.12]{LRU23}. 
Therefore, any nested ind-group 
is pseudo-divisible. 

5.
An element $g\in\Aut(X)$ is said to be 
\emph{algebraic}
if $g$ is contained in an algebraic 
subgroup of $\Aut(X)$. 
For an affine surface $X$, 
$g\in\Aut(X)$ 
is algebraic if and only if $g$ is 
pseudo-divisible,
see \cite[Corollary~2.6 
and~Theorem~3.1]{LRU23}. 
In 
\cite[Theorem~3.1]{LRU23} 
this equivalence is stated for the group 
$\Bir(X)$. 
The proof for $\Aut(X)$ goes as follows,
see the proof of 
Corollary~2.6 in [ibid]. 

Let $g\in\Aut(X)$.
For the implication 
\[g \text{ \emph {is algebraic}}
\Rightarrow g \text{ \emph{is 
pseudo-divisible}}\]
see the preceding remark. 
As for the converse, suppose $g$ 
is pseudo-divisible in $\Aut(X)$. Then 
it also is pseudo-divisible in $\Bir(X)$.
By \cite[Theorem~3.1]{LRU23}, 
the closure of 
$\langle g\rangle$ in  $\Bir(X)$
is an algebraic subgroup $G$ of $\Bir(X)$.
Since $G$ contains a dense subgroup 
$\langle g\rangle\subset\Aut(X)$,
by
\cite[Theorem~2]{Kra18}, $G$ 
acts regularly on $X$, that is, 
$G\subset\Aut(X)$ 
is an algebraic subgroup.
\erems
\blem\la{lem:alg-grp} 
For an affine algebraic group $G$ over $\kk$,
the following are equivalent:
\begin{itemize}
\item[$(i)$] $G$ is connected;
\item[$(ii)$] $G$ has no proper  
subgroup  
of finite index.
\end{itemize}
\elem
\bproof
If $G$ is disconnected, then $G^0$ is a 
proper subgroup of $G$ of finite index.
This proves the implication 
$(ii)\Rightarrow(i)$. 

For the converse, 
suppose that 
$G$ is connected and contains a 
proper subgroup $H$ of finite index. 
Every unipotent
$g\in G$
is contained in a $\Ga$-subgroup of $G$,
and every  semisimple $g\in G$ 
is contained in an algebraic subtorus of $G$,
see \cite[Theorem~22.2]{Hum75}. 
Consequently, these elements are divisible. 
It follows that all unipotent and 
semisimple elements of $G$
belong to~$H$. Since these 
elements generate 
$G$, we have $H=G$, 
which is a contradiction. 
\eproof 
By the following lemma, the implication 
$(i)\Rightarrow(ii)$ is also true for 
nested affine ind-groups. 
\blem \la{lem:divisible-1} 
Let $G=\varinjlim_i G_i$ be a nested 
affine ind-group. If $G$ is connected, then
it has no proper subgroup of finite index.
\elem 
\bproof 
Since $G$ is connected, for every  $g\in G$
there exists a closed affine curve $C\subset G$ 
such that $1, g\in C$. We have $C\subset G^0_i$
for some $i\in\NN$. It follows that 
$G=\varinjlim_i G^0_i$. 

Suppose there exists a proper 
subgroup of finite index $H\subset G$.
Then, for a sufficiently large index $i$,
$G^0_i\cap H$ is a proper finite-index 
subgroup of $G^0_i$.
This contradicts  Lemma \ref{lem:alg-grp}. 
\eproof
The example below shows 
that the equivalence
of Lemma \ref{lem:alg-grp} 
is not true in general for 
nested ind-groups.
 \bexa\la{ex:equiv} 
 Let $A_\infty=\varinjlim_i A_i$ 
 be the infinite alternating  group, 
 that is, 
 the group of all 
 even finite permutations of $\NN$. 
It is a countable nested affine ind-group.
Being a torsion group, $A_\infty$ 
is pseudo-divisible. 
Since the group $A_\infty$ is simple,
it has no proper finite-index subgroup.
However, it is disconnected.
 \eexa
Recall that the \emph{
$p$-primary component} of 
an abelian torsion group $G$
is the subgroup of all elements whose 
orders are powers of the prime number $p$. 
Every abelian torsion group is a 
direct sum of 
its $p$-primary components. 
We have the following criterion. 
 \blem\la{lem:tower} Let $C=\bigcup_i C_i$,
 where $C_1\subset C_2\subset \ldots$ 
 is a tower 
 of finite cyclic subgroups with 
 proper embeddings
 $C_i\subset C_{i+1}$. Then $C$ has 
 no proper subgroup of finite index 
 if and only 
 if every $p$-primary component 
 $C(p)$ of $C$ is infinite. 
 \elem
 \bproof  The $p$-primary component $C(p)$ 
 is finite
 if and only if the powers of $p$ that divide 
the orders $\ord(C_i)$ for some $i$ 
are bounded above.
 In this case, the sum of 
 all the primary components of $C$ 
that are not $C(p)$ is 
 a subgroup of finite index in $C$. 
 
To prove the converse, note that 
the $p$-primary component $C(p)$ 
is isomorphic to 
$\ZZ(p^\infty)=\varinjlim_k \ZZ/p^k\ZZ$, since
it is the union of an ascending chain 
of $p$-primary finite 
cyclic subgroups. It is well known that
 any proper subgroup of $\ZZ(p^\infty)$
 is finite. 
 
 Suppose all primary components $C(p)$ of $C$ 
 are infinite, and let $H\subset C$ be 
 a subgroup of 
 finite index. Then,  either $C(p)\cap H$ 
 is a finite group, 
 or $C(p)\supset H$. In the first case, $|C:H|$ 
 is infinite. 
If there is no such a component $C(p)$, 
then $H$ contains 
every $C(p)$, hence $H=C$.
 \eproof
 \bexas\la{ex:ab-torsion} $\,$
 
 1.
 Let $C=\varinjlim_k \ZZ/3\cdot 2^k\ZZ$. 
 Then the $2$-primary component
 $C(2)$ of $C$ has index $3$ 
 in $C\simeq\ZZ/3\ZZ\oplus C(2)$.
 
 2. Consider the abelian torsion 
 subgroup $\cR\subset\Gm$
 consisting of all roots of unity. We have 
 $\cR=\varinjlim_k \ZZ/k!\ZZ$. 
 Thus, $\cR$ is a nested ind-group. It
 is the direct sum over the primes $p\ge 2$
 of the $p$-primary components 
 $\cR(p)\simeq\ZZ(p^\infty)$. 
Obviously, $\cR\simeq\QQ/\ZZ$ is 
divisible and therefore 
has no proper subgroup of finite index. 
 \eexas
\section{Solvable connected  
subgroups of automorphisms}
\la{sec:solv-conn}
Let $X$ be an affine algebraic variety over $\kk$.
For the following proposition, 
see \cite[Proposition~3.10]{FP18}.
Since the argument is short, we provide it too. 
\bprop\la{prop:FP}
Let $G$ be a group endowed with a topology 
(and not necessarily a topological group). 
Suppose that the derived length of solvable 
subgroups of $G$ 
is bounded above.
Then every  solvable (resp.  solvable connected) 
subgroup of $G$ 
is contained in a maximal such subgroup.
\eprop
\bproof Since the derived lengths 
of the solvable subgroups 
of $G$ are bounded, for any chain $(G_i)_{i\in I}$ of
solvable  (connected) subgroups their 
union $\cup_i G_i$ is also solvable
(resp. connected).
We can therefore apply Zorn's lemma 
and conclude that any solvable  
(connected) subgroup 
is contained in a maximal such subgroup.
\eproof
\bdefi\la{def:Borel}
Let $G$ be an ind-group. 
A solvable connected subgroup $B$ of $G$ 
is called  a \emph{Borel subgroup} 
if $B$ is maximal under inclusion 
among the solvable 
connected subgroups of~$G$. 
\edefi
\brems\la{rem:Zariski} $\,$

1. 
The Zariski closure of a solvable 
(resp. connected) 
subgroup in 
an ind-group is again  a solvable 
(resp. connected) 
subgroup, see 
\cite[Lemma~ 2.3]{FK18} or 
\cite[Lemma~5.1.3(2)]{KZ24}. 
In particular, every Borel 
subgroup is closed. 

2. A Borel subgroup of a connected 
algebraic group $G$
is maximal among the solvable 
subgroups of $G$, 
see \cite[Corollary~23A]{Hum75}.

3. All Borel subgroups of 
an affine algebraic group are 
pairwise conjugate, 
see \cite[Theorem~21.3]{Hum75}. 

4. Since $\Aut(X)$ is an affine ind-group 
and any algebraic subgroup $G$ of $\Aut(X)$
is closed, $G$ is an affine algebraic group, 
see Remarks \ref{rem:connectedness}.1 and 2. 
\erems
Recall the following criterion for a
commutative connected  group
of automorphisms to be an algebraic group.
\bthm[{\rm\cite[Theorem~B]{CRX23}}]
\la{thm:CRX}
Let  $V$ be an irreducible 
algebraic subvariety of $\Aut(X)$ 
containing the identity. If every pair of elements
of $V$ commute, then the subgroup
$\langle V\rangle$ of $\Aut(X)$ generated 
by the elements of $V$ is a commutative 
algebraic subgroup of $\Aut(X)$.
\ethm
\bcor\la{cor:CRX} Let $G\subset\Aut(X)$ be
a connected commutative subgroup. 
\bnum
\item[$(a)$] If $G$ is a torsion group, 
then it is trivial. 
\item[$(b)$] Suppose that
$G\subset A*_{A\cap B} B$,  
where $A$ and $B$ 
are ind-subgroups of $\Aut(X)$.  
Let $T$ be the Bass--Serre tree 
of the latter amalgam. 
Then $G\subset\Aut(T)$ 
consists of elliptic elements.
\enum
\ecor
\bproof (a) Since $G$ is path connected, 
for every nonunit $f\in G$, there exists 
an irreducible affine 
algebraic curve $C$ in $G$ passing 
through $1$ and $f$. By  
Theorem \ref{thm:CRX},
$\langle C\rangle$ is 
a connected commutative affine 
algebraic group.  
Since all elements of $\langle C\rangle$ 
are  torsion elements, and 
therefore semisimple, 
$\langle C\rangle\simeq\G_{\mathrm m}^n$ 
is an algebraic torus. 
However, the multiplicative group 
$\G_{\mathrm m}$
contains elements of infinite order. 
Thus, $n=0$
and $\langle C\rangle=1$, 
which is impossible. 

(b) Let $f$, $C$ and $\langle C\rangle$ 
have the same 
meaning as before. By Serre's 
Theorem \ref{thm:Serre}, 
$\langle C\rangle$ is 
conjugate to a subgroup of either $A$ or $B$. 
So, $f$ is elliptic.  
\eproof
\bcor\la{cor:connect-solv} 
A solvable connected torsion subgroup 
$G\subset\Aut(X)$ is trivial. 
\ecor
\bproof Assume that $G\neq 1$.
The derived series 
\[G\supset G^{(1)}\supset\ldots
\supset G^{(n-1)}\supset 1\]
consists of connected  subgroups, 
where $G^{(n-1)}$
is a nontrivial commutative torsion subgroup. 
However, by Corollary \ref{cor:CRX}(a) 
it should be trivial, which is a contradiction. 
\eproof
\bdefi\la{def:tors} 
Let $G=A*_C B$ 
be a proper amalgam, where 
$C=A\cap B$. 
We say that a subgroup $H\subset G$
 is \emph{torsionally unbounded}
if every $f\in H$ acting on the Bass--Serre tree
$T:=T(G)$ with an unbounded fixed subtree $T^f$
is a torsion element. \edefi
\bthm \la{prop:clef}
Suppose that 
$G=A*_{C} B\subset\Aut(X)$,
and let 
$H\subset G$ be a
solvable connected  subgroup. If $H$ is 
 torsionally unbounded,
then $H$ is elliptic. 
\ethm
\bproof 
We begin with the following two claims. 

\medskip

\noindent {\bf Claim 1}. 
\emph{$H$ cannot be parabolic.}

\medskip

\noindent \emph{Proof of Claim 1}.
Suppose $H$ is parabolic. Then 
every $h\in H$ is elliptic and $H$ 
fixes a unique end
$b\in\partial T$,
see Proposition \ref{thm:Tits}(c).
For  $v\in T^h$, the
unique geodesic ray $r=[v,b[\in b$ with
a tip $v$
is contained in the  fixed subtree $T^h$.
Thus, $T^h$ is unbounded, and so
$h$ is a torsion element. 
Being a solvable connected  
torsion subgroup of $\Aut(X)$,
$H$ should be trivial by 
Corollary \ref{cor:connect-solv}, 
which is a contradiction.
\qed

\medskip

\noindent {\bf Claim 2}. 
\emph{Suppose $H$ is loxodromic. Then 
every elliptic $f\in H$ is a torsion element.}

\medskip

\noindent \emph{Proof of Claim 2}.
Let
$f\in H$ be elliptic. Since $H$ is solvable, 
$f$ either fixes a point of $\partial T$, or leaves 
invariant a pair of points of $\partial T$,
see Proposition \ref{thm:Tits}(b). 
Anyway, $f^2$ 
fixes a point $b\in\partial T$. For $v\in T^f$,
$f^2$ fixes $v$ and $b$, and so it
fixes the geodesic ray $r=[v,b[$ pointwise. 
Hence $T^{f^2}$ is unbounded, and therefore
$f^2$ and $f$ are torsion elements. 
\qed

\medskip

Returning to the proof of the theorem,
suppose, on the contrary, that $H$ is not elliptic. 
Then, 
by Claim 1, $H$ is loxodromic. 
Due to Claim 2, every elliptic $f\in H$
is a torsion element. 

Let  $H^{(n-1)}$ be the penultimate member 
of the derived series of $H$. 
Then, $H^{(n-1)}$ is 
a nontrivial connected commutative subgroup
consisting of elliptic elements, 
see Corollary \ref{cor:CRX}(b).
Since this is a torsion subgroup, 
it should be trivial by 
Corollary \ref{cor:CRX}(a). 
This contradiction 
ends the proof. 
\eproof
The following corollary is immediate. 
\bcor\la{cor:clef} Suppose that 
$G=A*_{A\cap B} B\subset\Aut(X)$
 is torsionally unbounded.
Then every Borel subgroup of $G$ is 
conjugate to a Borel subgroup of $A$ 
or $B$. 
\ecor
The next lemma generalizes 
\cite[Lemma~3.7]{Lam01}, see also
\cite[Lemma~7.55]{Lam25}.
\blem\la{lem:commut-ell}
Suppose that 
$G=A*_{A\cap B} B\subset\Aut(X)$
 is torsionally unbounded and $C=A\cap B$ is 
 a connected solvable algebraic group. 
 Then every two elements
 $f_1,f_2\in G$ with unbounded 
 subtree $T^{f_1}\cap T^{f_2}$ 
 commute. In particular, 
any parabolic subgroup $P\subset G$
 is an abelian torsion group.
 The same is true for
 the subgroup $F$ of 
 an elementary loxodromic 
 subgroup $H\subset G$, where $F$ consists
 of all elliptic elements of $H$
which  act identically on the axis of $H$.  
\elem
\bproof It suffices to prove the first assertion; 
 then the other two easily follow. 

The commutator $h=[f_1,f_2]$ 
acts identically on $T^{f_1}\cap T^{f_2}$,
hence it is a torsion element. Since the $G$-action 
on the set of edges of $T$ is transitive, 
up to conjugation we can assume that
the edge $1\cdot C$ is fixed by $f_1$ 
and $f_2$, hence also by $h$. Thus,
we have $f_1,f_2\in C$ and $h\in [C,C]=U$,
where $U$ is the unipotent radical of $C$,
see \cite[Sections~17.3 and~17.6]{Hum75}. 
It follows that $h=1$.
\eproof
In Section \ref{sec:toric}
we use the following lemma.
\blem\la{lem:quotient} Let $G=A*_{C} B$,
 where $C=A\cap B$,
and let $Z\subset C$ be 
a central torsion subgroup of $G$. 
Let $\pi\colon G\to G':=G/Z$ 
be the canonical surjection.
Consider the amalgam 
$G'=A'*_{C'} B'$, where
$A'=\pi(A)$, $B'=\pi(B)$, 
and $C'=\pi(C)$. 
Then $G$ is torsionally 
unbounded if and only if
 $G'$ is.
\elem
\bproof Let $T'=T(G')$ 
be the Bass-Serre tree of 
the amalgam $G'=A'*_{C'} B'$.
Then $\pi$ induces a map 
of graphs $\pi_*\colon T\to T'$, where
$\pi_*\colon gA\mapsto \pi(g)A'$, 
$gB\mapsto \pi(g)B'$, and 
$gC\mapsto \pi(g)C'$. 
It is easy to see that 
\begin{itemize}
\item $\pi_*$ is an isomorphism of trees, 
equivariant with respect to the $G$-action
on $T$ and the $G'$-action on $T'$;
\item 
a vertex 
$v=gA\in \vert(T)$
is fixed by $f\in G$ if and only if 
the vertex $\pi_*(v)=\pi(g)A'\in \vert(T')$ 
is fixed by $\pi(f)$;
\item
$f$ and $\pi(f)$ 
 are simultanuously torsion elements, or 
non-torsion elements. 
\end{itemize}
Now our assertion follows.
\eproof
The following corollary is 
a direct consequence of the proof. 
\bcor\la{cor:quotient} A subgroup 
$H'\subset G'$ is elliptic resp. 
parabolic resp. 
(elementary, focal, general) loxodromic 
if and only if $H=\pi^{-1}(H')$ is. 
\ecor
\section{$\Aut(\AA^2)$ and its 
subgroups}\la{sec:A2}
In this section we collect 
some results about subgroups of 
$\Aut(\AA^2)$ that will be useful 
in the next sections. See e.g.
\cite{BEE16}, \cite{BS15}, 
\cite[Section~16]{FK18}, \cite{FL10},
\cite{FP18}, \cite{Lam01},
\cite{Lam25}, and \cite{LRU23} 
for additional information. 
\subsection{Generalities}
\bnota
We denote by $\JONQ^+(\AA^2)$ 
($\JONQ^-(\AA^2)$,
respectively) the group 
of de Jonqi\`eres
transformations
 \be\label{J-plus}
 \Phi^+:(x,y)\longmapsto (\alpha
x+p(y),\beta y+\gamma),
\ee 
respectively,
\be\label{J-minus}\Phi^-:(x,y)\longmapsto 
(\alpha x+\gamma,\beta
y+p(x)),
\ee 
where 
$\alpha,\beta\in\mathbb{G}_{\mathrm m}$, 
$\gamma\in\kk$,
and $p\in\kk[t]$. The subgroup 
\be\la{aff-pm}
{\Aff}^{\pm}(\AA^2)=
\Aff(\AA^2)\cap{\JONQ}^{\pm}(\AA^2)
\ee 
consists of all
upper (respectively lower) triangular 
affine transformations
$\Phi^\pm$ with $\deg p\le 1$.
We let 
\[B^\pm={\Aff}^{\pm}(\AA^2)\cap\GL(2,\kk).\]
\enota
\brems\la{rem:deg} $\,$ 

1. 
According to Jung--van der 
Kulk--Nagata's theorem, 
we have
\be\label{eq:JvdK}
\Aut(\AA^2)={\JONQ}^+(\AA^2)*
_{{\Aff}^+(\AA^2)}
\Aff(\AA^2)=\langle{\JONQ}^+(\AA^2)
,\,\tau\rangle,
\ee
where
 \[\tau\in\Aff(\AA^2), \quad
\tau:(x,y)\longmapsto (y,x)\] 
is a twist. 
Note that $\tau\in\langle B^+,\,
B^-\rangle$ and
\be\la{eq:plus-minus}
{\JONQ}^-(\AA^2)
=\tau{\JONQ}^+(\AA^2)\tau,
\quad {\Aff}^-(\AA^2)
=\tau{\Aff}^+(\AA^2)\tau.\ee
In particular,
\[\Aut(\AA^2)=\langle{\JONQ}^+(\AA^2),\,
{\JONQ}^-(\AA^2)\rangle
=\langle{\JONQ}^+(\AA^2),\,\tau\rangle
=\langle U^+, \Aff(\AA^2)\rangle,\]
where 
\be\la{max-to} U^+=
\{\Phi\in{\JONQ}^{+}(\AA^2)\,|\,
\Phi:(x,y)\mapsto
(x+p(y),y),\quad p\in\kk[y]\}
\ee 
is the unipotent radical 
of ${\JONQ}^+(\AA^2)$ and $U^-=\tau
U^+\tau$ is that of ${\JONQ}^-(\AA^2)$.
However, 
$\Aut(\AA^2)\neq U^+ *_{
U^+\cap\Aff(\AA^2)} \Aff(\AA^2)$, 
see Remark in \cite[\S 2]{Kam75}.

2. For an endomorphism $\varphi=
(p,q)\in \End(\AA^2)$, 
where $p,q\in\kk[x,y]$, we set 
\[\deg(\varphi)=
\max\{\deg(p),\deg(q)\}.\]
If $\varphi\in\Aut(\AA^2)
\setminus\Aff(\AA^2)$ 
is written as 
\[\varphi=a_1\tau a_2\tau\cdots \tau  
a_n,\quad n\ge 1,\]
where
$a_i
\in{\JONQ}^+(\AA^2)$, 
cf. \eqref{eq:JvdK} 
and \eqref{eq:plus-minus},
then
\be\deg(\varphi)=\deg(\varphi^{-1})=
\prod_{i=1}^n \deg(a_i).\ee

3. 
The standard ind-structure on $\End(\AA^2)$ is 
defined by ascending filtration
\[\End(\AA^2)=\bigcup_{d=0}^\infty 
\End(\AA^2)_d,
\quad\text{where}\quad 
\End(\AA^2)_d
=\{\varphi\in\End(\AA^2)\,|\,
\deg(\varphi)\le d\}\]
are affine algebraic varieties 
and the inclusions
$ \End(\AA^2)_d\subset 
\End(\AA^2)_{d+1}$ 
are closed embeddings. 
The cartesian square $\End(\AA^2)^2$
inherits an ind-structure, 
along with its closed subvariety
\[\cV=\{(\varphi,\psi)\in 
\End(\AA^2)^2\,|\,\varphi\psi=1\}.\]
An ind-structure on $ \Aut(\AA^2)$ 
can be defined by a
bijection 
\[ \Aut(\AA^2)\ni\varphi\mapsto 
(\varphi,\varphi^{-1})\in \cV,\]
see \cite[Theorem~5.1.1]{FK18}.

\smallskip

4. The group
${\JONQ}^{\pm}(\AA^2)=
\bigcup_d {\JONQ}^{\pm}(\AA^2)_d$ 
is a solvable nested ind-subgroup 
of $\Aut(\AA^2)$
filtered by the solvable connected  
algebraic subgroups 
\[{\JONQ}^{\pm}(\AA^2)_d
:=\{\Phi^\pm\in {\JONQ}^{\pm}(\AA^2)\,|\,
\deg(p)\le d\}.\]
\erems
By Serre's Theorem \ref{thm:Serre}, 
see also \cite[Theorem~4.3]{Kam79} 
and \cite[Proposition~1.11]{Wri79},
we have:
\bthm\label{Serre}
Any subgroup of $\Aut(\AA^2)$ 
of bounded length (in particular, 
any algebraic subgroup)
is conjugate to a subgroup 
of one of the factors 
${\JONQ}^+(\AA^2)$ and $\Aff(\AA^2)$  
in \eqref{eq:JvdK}.
\ethm
\bnota\la{not:Gde}
Let $\zeta\in\Gm$ be  
a primitive root of unity 
of order $d$, and let $e\in\NN$ verifies
$1\le e< d$ and $\gcd(e,d)=1$. Consider
the cyclic subgroup $G_{d,e}\subset\TT$
generated by 
$f\colon (x,y)\mapsto (\zeta^{e} x,\zeta y)$.
Clearly, any 
$\phi\colon (x,y)\mapsto (\alpha x, \beta y)$,
where $\alpha$ and $\beta$ 
are primitive roots of unity of order $d$,
belongs to some $G_{d,e}$, 
and any subgroup $H$ of $G_{d,e}$ 
is of the form $H=G_{d',e'}$, 
where $d'$ divides $d$ and 
$e'\equiv e\mod d'$. 
\enota
For the action of $\Aut(\AA^2)=A*_C B$
on its Bass-Serre tree $T$, 
we have the following results.
\bprop[{\rm \cite{Lam01}, 
\cite[Chapter~7]{Lam25}}]
\la{prop:Lamy} Let $f,f_1,f_2\in \Aut(\AA^2)$ 
be elliptic and $g,g_1,g_2\in \Aut(\AA^2)$ be
loxodromic.
$\,$

\bnum
\item[$(a)$] For any $f\in G_{d,e}$, 
there is a loxodromic $g\in\Aut(\AA^2)$ 
that commutes with $f$.
\item[$(b)$] 
The fixed subtree $T^f$ is unbounded
if and only if $f$ is conjugate to 
an element of some $G_{d,e}$. 
In particular, $\Aut(\AA^2)$ is 
torsionally unbounded. 
 \item[$(c)$] 
If $T^{f_1}\cap T^{f_2}$ is unbounded, then
$f_1$ and $f_2$ commute, 
and $f_i$ fixes a unique point $P_i\in\AA^2$, 
where $P_1=P_2$. 
\item[$(d)$]  If $T^f\cap\Ax(g)$ is unbounded,
then $f$ commutes with $g^2$ 
and $T^f\supset\Ax(g)$.
\item[$(e)$]  If $\Ax(g_1)\cap \Ax(g_2)$
is unbounded, then 
$\Ax(g_1)=\Ax(g_2)$ and $g_1$ 
and $g_2$ commute. 
\item[$(f)$]  If $f$ and $g$ commute, 
then $\Ax(g)\subset T^f$.
\enum
\eprop
\bproof See \cite[Proof of Lemma~7.50(2)]{Lam25}
for (a), \cite[Proposition 7.54]{Lam25} for (b), 
\cite[Lemma~7.55]{Lam25} for (c), 
\cite[Lemma 7.56]{Lam25} for (d), 
\cite[Lemma 7.57]{Lam25} for (e),
and Proposition~\ref{prop:elements}(b) for~(f).
\eproof
\subsection{Classification of subgroups 
of $\Aut(\AA^2)$}
In this subsection we provide Lamy's 
classification of subgroups of 
$\Aut(\AA^2)$. We let 
$T=T(\Aut(\AA^2))$ be
the Bass-Serre tree 
that corresponds to the amalgam 
in \eqref{eq:JvdK}. 
\bthm[{\rm\cite[Th\'eor\`eme~2.4, 
Propositions 3.12 and~4.10]{Lam01}, 
\cite[Theorem~7.60 
and Proposition~7.61]{Lam25}}]
\la{thm:Lamy0}
 Consider a subgroup 
$H\subset\Aut(\AA^2)$ acting on $T$. 
\bnum
\item[$(a)$] 
If $H$ is parabolic, 
then $H$ is a union of 
 a strictly increasing sequence 
 of finite cyclic subgroups $C_n$, 
 where $C_n$ is conjugate to some
 $G_{d_n,e_n}$. 
 Consequently, 
 $H$ is a countable abelian torsion
 group, it is not finitely generated and fixes 
 a unique end $b\in\partial T$. 
\item[$(b)$]
If $H$ is  elementary loxodromic, 
then $H$ is solvable and contains
a subgroup $H_0\subset H$ 
of index at most 2
of the form 
$H_0=F\rtimes\langle g\rangle\simeq 
F\rtimes\ZZ$, 
where $g\in H_0$ is loxodromic and 
$F\subset H_0$ 
is a finite abelian subgroup.
\item[$(c)$]  $H$ cannot be focal loxodromic. 
\item[$(d)$] If $H$ is  general loxodromic, 
then $H$ contains a free subgroup 
of rank $2$.
\enum
\ethm
\brem Let $H\subset\Aut(\AA^2)$ 
be elementary loxodromic
and $F$, $H_0$ and $g$ be as in (b).
Then any $f\in F$ fixes ${\rm Ax}(g)$ 
pointwise, and 
therefore, it is conjugate to 
an element of some~$G_{d,e}$,
see Proposition \ref{prop:Lamy}(b).
If $H\neq H_0$, then any 
$f\in H\setminus H_0$ acts on ${\rm Ax}(g)$ 
by an involution fixing a vertex 
of ${\rm Ax}(g)$.
\erem
\subsection{Borel subgroups 
of $\Aut(\AA^2)$}
In the following theorem, we recall
some important properties of solvable 
resp. Borel subgroups 
of $\Aut(\AA^2)$. These
properties will be used
later
 to establish similar results for the
subgroups of $\Aut_0(\AA^2)$ and 
$\Aut(X_{d,e})$.
\bthm\la{thm:Borel-A2}
\label{thm:BEE} $\,$
\bnum\item[$(a)$] 
$\JONQ^+(\AA^2)$ is a Borel 
subgroup of $\Aut(\AA^2)$. It is 
solvable of derived length $3$.
\item[$(b)$] 
Every Borel subgroup of $\Aut(\AA^2)$ 
is conjugate to $\JONQ^+(\AA^2)$. 
\item[$(c)$] Every solvable 
connected  subgroup 
of $\Aut(\AA^2)$ 
is contained in a Borel subgroup. 
In particular, 
it is conjugate to a subgroup of 
$\JONQ^+(\AA^2)$.
\item[$(d)$] Every solvable subgroup 
of $\Aut(\AA^2)$ has derived length $\le 5$ and
is contained in a maximal solvable subgroup. 
The upper bound $5$ is optimal.
\item[$(e)$] $\JONQ^+(\AA^2)$ is maximal 
among solvable subgroups 
of $\Aut(\AA^2)$ and among 
proper closed subgroups of $\Aut(\AA^2)$. 
It coincides with its normalizer 
in $\Aut(\AA^2)$.
\item[$(f)$] 
A non-abelian maximal solvable subgroup 
$H$ of $\Aut(\AA^2)$
is a Borel subgroup if and only if 
$H$ contains no proper subgroup 
of finite index. 
\item[$(g)$] 
Any subgroup $G\subset\Aut(\AA^2)$
that fixes an end $b\in\partial T$ is solvable 
of derived length at most $3$. 
\enum
\ethm
\bproof The analogues of (a), (b), and (f)  
for the subgroup 
$G=\{f\in\Aut(\AA^2)|{\rm Jac}(f)=1\}$
are proved in 
\cite[Theorems~1 and~4]{BEE16}.
The proofs for $\Aut(\AA^2)$ are similar;
cf. \cite[Lemma~3.2]{FP18} for (a). 
See 
\cite[Propositions~3.4,~3.5,~3.10,~3.14,~3.16,
Theorem~3.8 and Corollary~3.15]{FP18}
for (c), (d) and (e); 
cf. also \cite[Proposition~5(e)]{BEE16}.
The proofs in \cite{FP18} 
are done for 
$\kk=\CC$, but they apply mutatis mutandis
in our more general setup.
An alternative proof of the first 
statement of (e) can be found  in
 \cite[Corollary~7.62]{Lam25}. 
 Yet another proof can also 
be deduced from 
 Proposition \ref{thm:maximal}. 
 See \cite[Lemma~7.64]{Lam25} for (g). 
 \eproof
 \brem The closedness condition 
 in (e) is important. 
 Indeed, $\JONQ^+(\AA^2)$
  is not maximal among proper subgroups of 
  $\Aut(\AA^2)$, see \cite[Theorem~1.4]{FP18}.
 \erem
\subsection{Maximal solvable 
subgroups of $\Aut(\AA^2)$}
\la{ss:max-solv-A2} 
\subsubsection{Elliptic subgroups}
Let us start with the following remarks.
\brems\la{rem:exas} $\,$

1. By definition, an elliptic solvable 
subgroup either is 
contained in a Borel subgroup, that is,
is conjugate to a subgroup 
of $\JONQ(\AA^2)$,
or is conjugate to a solvable 
subgroup of $\Aff(\AA^2)$.

2. 
Any solvable subgroup $S$ of derived length 
$\ge 4$ in $\Aut(\AA^2)$ is conjugate 
to a subgroup of $\Aff(\AA^2)$. 
Indeed, the parabolic subgroups
are abelian, see Theorem \ref{thm:Lamy0}(a).
Any elementary loxodromic subgroup
has derived length $\le 3$, 
see \cite[Proof of Proposition~3.14]{FP18} or
Proposition \ref{rem:case-A2}(a) below.
The general loxodromic 
subgroups of $\Aut(\AA^2)$ 
are not solvable, see Theorem~\ref{thm:Lamy0}(d)
or Proposition \ref{thm:Tits}(e).
Therefore, $S$ is elliptic. 
It cannot be conjugate to 
a subgroup of $\JONQ^+(\AA^2)$,
because  
the derived length of $\JONQ^+(\AA^2)$ 
is equal to 3. 
Here is our assertion. 

3. 
There is a solvable subgroup
$S$ of $\Aff(\AA^2)$ 
with derived length $5$, see 
\cite[Section~3.2]{FP18}; its derived group 
$S^{(1)}\subset\GL(2\kk)$ 
has derived length $4$.
By  the above argument, 
these subgroups are not contained 
in any Borel subgroup.

4. In Proposition \ref{prop:max} 
we give a list of all elliptic maximal 
solvable subgroups of $\Aut(\AA^2)$. 
This list includes an
algebraic subgroup 
of $\Aff(\AA^2)$ 
with derived length $3$
that is not contained 
in any Borel subgroup of 
$\Aut(\AA^2)$. 
It is metabelian and 
contains a propers subgroup 
of index~2.
\erems
\subsubsection{Parabolic and 
loxodromic subgroups}
In the following proposition  
we describe the parabolic and loxodromic 
maximal solvable subgroups 
of $\Aut(\AA^2)$. Being disconnected 
 and maximal solvable, these subgroups
are not contained in 
 any Borel subgroup of $\Aut(\AA^2)$.
\bprop\la{rem:case-A2} $\,$
\bnum
\item[$(a)$]
Given a loxodromic $g\in \Aut(\AA^2)$, 
consider 
the elementary
 loxodromic subgroup $H(g)$
 consisting of all 
$h\in \Aut(\AA^2)$ that leave invariant 
$\Ax(g)$. 
Then $H(g)$ is a maximal solvable subgroup 
of $\Aut(\AA^2)$ that
has derived length $\le 3$, and 
any loxodromic
subgroup $H\subset\Aut(\AA^2)$ 
which is maximal solvable
coincides with $H(g)$
for some $g\in H$. 
\item[$(b)$] Let $P\subset\Aut(\AA^2)$ 
be a parabolic subgroup, 
$b\in\partial T$ 
be the unique  end  fixed by $P$, 
and $P_b\supset P$
 be the parabolic subgroup consisting 
of all elliptic $f\in \Aut(\AA^2)$ 
fixing $b$. 
Then $P_b$ is a maximal 
solvable subgroup 
 of  $\Aut(\AA^2)$ of derived length $1$. 
 Furthermore, every parabolic  
 subgroup of
 $\Aut(\AA^2)$
that is maximal solvable 
 coincides with $P_b$
 for some $b\in\partial T$.
\enum
\eprop
\bproof (a) By  Lamy's Theorem 
\ref{thm:Lamy0}(b), 
 $H$ contains a subgroup 
 $H_0$ of index $\le 2$
which is 
 an extension of $\ZZ$ by 
 a  finite abelian subgroup $F$.
 Thus, the derived length of $H$ is 
 $\le 3$. 
 
 Since 
$\Aut(\AA^2)$ contains no focal 
loxodromic subgroup, see Theorem 
\ref{thm:Lamy0}(c), $H=H(g)$ 
  is a maximal 
  solvable subgroup of $\Aut(\AA^2)$ 
  by Corollary \ref{cor:lox-solv2}.
 
 A loxodromic maximal solvable subgroup 
 $H\subset\Aut(\AA^2)$ cannot be focal 
 loxodromic nor general loxodromic, 
 see Theorem \ref{thm:Lamy0}(c) and (d).
 Therefore, $H$ is elementary loxodromic,
 and so $H=H(g)$ for 
 a loxodromic element $g\in H$.

Statement (b) follows from
 Corollary \ref{cor:lox-solv2} 
 and Theorem \ref{thm:Lamy0}(a).
\eproof
\brems\la{rem:Mold-BEE} $\,$

1. The group $\Aut(\AA^2)$ contains 
plenty of loxodromic maximal solvable 
subgroups. Indeed, by (a), any element  
of $\Aut(\AA^2)$ of positive even length
is contained in such a subgroup $H(g)$.
 
2.
Parabolic abelian
subgroups first appeared
in a paper by Moldavanskii, who
classified abelian subgroups 
of an amalgam, 
see \cite[Theorem~2]{Mol67}. 
Wright, in \cite{Wri79} interpreted 
Moldavanskii's
classification in terms of 
the Bass-Serre theory.
In the case of
Jung--van der Kulk--Nagata's 
amalgam structure on $\Aut(\AA^2)$,
he presented in \cite[Section~2]{Wri79} 
examples of parabolic  abelian 
subgroups of $\Aut(\AA^2)$ 
acting on the Bass-Serre tree 
$T=T(\Aut(\AA^2))$, and among them, 
a parabolic subgroup $P$ isomorphic to the group
$\cR\subset\Gm$ of all roots of unity, see
\cite[Theorem~2.4]{Wri79}). See also
 \cite[Example~7.66]{Lam25} for Wright's examples. 
 
3.
 In a remark  following \cite[Theorem~4]{BEE16},
 the authors
notify the existence of an abelian subgroup 
 of $\Aut(\AA^2)$ that is maximal solvable 
 and contains no  finite-index
 subgroups~\footnote{We 
 are grateful to A.~Eshmatov, 
 who kindly provided us with 
 additional information concerning 
 this subgroup.}.
These examples force to impose
 the assumption of non-abelianess in 
 \cite[Theorem~4]{BEE16};
see also Theorem \ref{thm:Borel-A2}(f).
Indeed,
 Wright's abelian parabolic subgroup $P$ 
verifies the assumptions of this theorem, 
see the next lemma. 
However, being parabolic, 
it is not contained in a Borel subgroup.
\erems
\blem\la{exa:Wright-BEE}
Let $P\subset\Aut(\AA^2)$ be
a parabolic subgroup 
isomorphic to $\cR\cong\QQ/\ZZ$. 
Then $P$ is a maximal solvable subgroup 
that contains no proper
subgroup of finite index 
and is not a
Borel subgroup of $\Aut(\AA^2)$.
\elem
\bproof
Since $P\simeq\cR$ is divisible, it has no 
proper subgroup of finte index,
see Lemma \ref{lem:divisible}(b) 
and Example 
\ref{ex:ab-torsion}.2.
Since $P$ is parabolic, 
it fixes a unique end $b\in\partial T$,
see Proposition~\ref{thm:Tits}(c).
Let $P_b\supset P$ be 
the parabolic subgroup consisting of 
all elliptic $f\in\Aut(\AA^2)$ fixing $b$.
By Proposition \ref{rem:case-A2}(b), 
$P_b$ is maximal 
among solvable subgroups of $\Aut(\AA^2)$.
We assert that $P_b=P$. 

Indeed, when decomposing $P\simeq\QQ/\ZZ$ 
into a direct sum of $p$-primary components, 
all terms are of the form $\ZZ(p^\infty)$,
and for every prime number $p$, 
there exists exactly 
one term of this type, cf. Example 
\ref{ex:ab-torsion}.2.
Since $P_b$ contains $P$, the 
$p$-primary component 
$P_b(p)$ contains $P(p)\simeq\ZZ(p^\infty)$.
Recall that $P_b$ is the union of 
an ascending chain 
of finite cyclic groups, 
see Theorem \ref{thm:Lamy0}(a).
Consequently, every primary component 
$P_b(p)$ of $P_b$ is also such a union.
It follows that $P_b(p)\simeq\ZZ(p^\infty)$ 
for any  prime number $p$. 
Therefore, $P_b=P$. 
Since $P$ is a
maximal solvable  torsion subgroup,
it is not conjugate to $\JONQ^+(\AA^2)$. 
In particular, $P$ is not 
a Borel subgroup  of $\Aut(\AA^2)$.
\eproof
\subsection{$\Aut_0(\AA^2)$ and its subgroups}
Let 
\[{\Aut}_0(\AA^2)=\{f\in\Aut(\AA^2)\,|\,f(0)=0\}.\]
The following analogue
of the Jung--van der Kulk--Nagata 
theorem holds, 
see e.g.  \cite[Lemma~7.30]{Lam25}.
\blem\la{lem:amalg-A0}
There is an induced amalgam
(see Definition \ref{induced-amalg})
$\Aut_0(\AA^2)=A_0*_{C_0} B_0$ with
\be\la{eq:decomp-0} A_0=\GL(2,\kk),\quad 
B_0={\Jonq}^+(\AA^2),
\quad\text{and}\quad 
C_0=A_0\cap B_0=B^+,\ee
where $B^+\subset \GL(2,\kk)$ stands for 
the Borel subgroup 
of upper triangular matrices and the
$\Jonq^\pm(\AA^2):= \JONQ^\pm(\AA^2)
\cap{\Aut}_0(\AA^2)$ 
are consisting,
respectively, of the transformations 
\be\la{eq-Jonq-plus}
\varphi^+:(x,y)\longmapsto 
(\alpha x+yp(y),\beta y)
\ee 
and
\be\la{eq-Jonq-minus} 
\varphi^-:(x,y)\longmapsto (\alpha x,\beta
y+xp(x)) 
\ee 
with $\alpha,\beta\in\mathbb{G}_{\mathrm m}$ 
and $p\in\kk[t]$. \footnote{
Therefore $\Phi^\pm$ as in (\ref{J-plus}), 
(\ref{J-minus}) belongs to
$\Jonq^\pm(\AA^2)$ 
if and only if $\gamma=0$ and $p(0)=0$.}
\elem
We have a decomposition
\be\la{eq-unip} 
{\Jonq}^{+}(\AA^2)\simeq
U_0^+\rtimes\TT
\ee 
where
\[ \TT=\{\delta\in{\Jonq}^+(\AA^2)\,|
\,\delta\colon (x,y)\mapsto
(\alpha x, \beta y)\}\subseteq\GL(2,\kk)\]
  is the maximal torus  
  of ${\Jonq}^{+}(\AA^2)$
and the infinite dimensional commutative vector
group 
\[U_0^+=\{(x,y)\mapsto 
(x+yp(y),y)\}\simeq (\kk[t],+)\]
is the 
unipotent radical  of ${\Jonq}^{+}(\AA^2)$. 
There is a canonical surjection 
$\rho:{\Jonq}^{+}(\AA^2)\to\TT=
{\Jonq}^{+}(\AA^2)/U_0^+$.
\blem[{\rm \cite[Lemma~2.1]{AZ13}}]
\la{lem-uni-rad}$\,$
\bnum\item[$(a)$]  
$\varphi^+\in{\Jonq}^{+}(\AA^2)$  
as in $(\ref{eq-Jonq-plus})$ 
is semi-simple if and only if 
there
exists $u\in U_0^+$ such that $u^{-1}
\varphi^+u=\delta\in\TT$,
where $\delta=\rho(\varphi^+)$.
\item[$(b)$] $\varphi^+
\in {\Jonq}^{+}(\AA^2)$ is of
finite order if and only if it is semi-simple and
$\delta:=\rho(\varphi^+)\in\TT$ is of finite order.
 \enum
  \elem
  \blem[{\rm \cite[Proposition~2.5]{AZ13}}]
\la{prop-quasitorus} 
Any reductive algebraic subgroup 
$G\subset {\Jonq}^{+}(\AA^2)$ is
conjugate to a subgroup of the torus $\TT$. 
More precisely, there
exists $u\in U_0^+$ such that 
 $G^{u}= \rho(G)\subseteq\TT$ 
 where $G^{u}=u Gu^{-1}$.
 \elem
  \blem[{\rm \cite[Lemma~2.4]{AZ13}}]
  \la{lem-commuting} Let
$g,\tilde{g}\in {\Jonq}^{+}(\AA^2)$ 
be given by
\[g\colon (x,y)\mapsto (\alpha x+
p(y),\beta y)\quad\text{and}\quad
\tilde{g}\colon (x,y)\mapsto (\tilde{\alpha}
x+\tilde{p}(y),\tilde{\beta} y)\] 
with
\[p(y)=\sum_{m\ge 1}a_my^m\quad\text
{and}\quad \tilde{p}(y)=\sum_{m\ge
1}\tilde{a}_my^m.\] 
Then, $g$ and $\tilde{g}$ commute
if and only if 
\[\la{eq-commuting}
a_m({\tilde{\beta}}^{m}-\tilde{\alpha})=
\tilde{a}_m(\beta^{m}-\alpha)\quad\forall
m\ge 0.\]
\elem
The following corollary is immediate. 
\bcor\la{cor:commut} 
Let $f\colon (x,y)\mapsto (\zeta^e x, \zeta y)$, 
where $\zeta\in\kk$ 
is a primitive root of unity of order $d$ 
and $1\le e<d$, $\gcd(d,e)=1$, and let 
$1\le e'<d$, $ee'\equiv 1\mod d$.
Then, for the centralizer of $f$ in 
${\Jonq}^{\pm}(\AA^2)$ we have
\[{\Cent}_{{\Jonq}^{+}(\AA^2)}(f)=
\{(x,y)\mapsto (\alpha
x + y^{e}p(y^d), \, \beta y) \,|\, p\in\kk[t], 
\,\,\,\alpha,\beta\in\mathbb{G}_{\mathrm m}\}\]
and
\[{\Cent}_{{\Jonq}^{-}(\AA^2)}(f)=\{(x,y)
\mapsto (\alpha
x, \, \beta y + x^{e'}q(x^d)) \,|\, q\in\kk[t], 
\,\,\,\alpha, \beta\in\mathbb{G}_{\mathrm m}\}.\]
\ecor 
Let $T_0$ be the Bass--Serre tree 
of $\Aut_0(\AA^2)$.
Being an amalgam 
(see Lemma \ref{lem:amalg-A0}), 
the group $\Aut_0(\AA^2)$
 acts without inversions on $T_0$. 
 This action is faithful, see 
 \cite[Remarque~3.5]{Lam01}.
\bcor[{\rm\cite[Lemma~7.50(2)]{Lam25}}] 
\la{cor:ell-fix} Let $f$ be as in Corollary 
\ref{cor:commut}. Then
there exists a loxodromic 
$g\in\Cent_{\Aut_0(\AA^2)}(f)$. 
Consequently, the fixed point 
set $T_0^f\supset {\rm Ax}(g)$ 
is unbounded. 
\ecor
For the reader's convenience,
we provide an argument.
\bproof The second assertion 
follows from the first 
because of Proposition \ref{prop:elements}(b).
To prove the first, consider the elements
\[\tau\in A_0\setminus B_0= \GL(2,\kk)\setminus 
B^+\quad\text{and}   
\quad g_n^+\in B_0\setminus A_0=
{\Jonq}^+(\AA^2)\setminus 
B^+,n\ge 2,\]  
where
\[\tau \colon (x,y)\mapsto (y,x)
\quad\text{and}\quad
g_n^+\colon (x,y)\mapsto (x+y^n,y).\]
If $e=1$, then we take 
\[g=\tau g_{d+1}^+\colon 
(x,y)\mapsto (y,x+y^{d+1}).\]
If $e>1$, then we take 
$g=g_{e'} g_{e}=
\tau g_{e'}^+\tau g_{e}^+$, 
where
\[g_{e'}=\tau g_{e'}^+\colon 
(x,y)\mapsto (y,x+y^{e'}) \quad 
\text{and}\quad g_{e}=\tau g_{e}^+\colon (x,y)
\mapsto (y,x+y^{e}).\] 
In both cases $g\in\Aut_0(\AA^2)$ 
of even length
is loxodromic and commutes with $f$.
\eproof
\bprop\la{prop:unbound}
Let $f \in \Aut_0(\AA^2)$ be elliptic. 
Then $T_0^f$ is unbounded if and only if 
$f$ is conjugate to 
$(x,y)\mapsto (\zeta^e x,\zeta y)$, where 
$\zeta$  is a root of unity of degree $d>1$ 
and $e\in\NN$, $1\le e< d$ 
and $\gcd(e,d)=1$. 
Thus, $\Aut_0(\AA^2)$ is 
torsionally unbounded. 
\eprop
\bproof
An analogous statement is true
for the Bass--Serre tree $T$
of $\Aut(\AA^2)$, see  
\cite[Proposition~7.54]{Lam25}. 
For $\Aut_0(\AA^2)$ acting 
on $T_0\subset T$ 
 (see Remark \ref{rem:tree-embed})
our assertion follows,
because the loxodromic elements 
as in the proof of Corollary \ref{cor:ell-fix} 
fix the origin of $\AA^2$.
\eproof
Lamy's classification of subgroups of 
$\Aut(\AA^2)$ (see Theorem \ref{thm:Lamy0})
applies mutatis to $\Aut_0(\AA^2)$.

Due to the following corollary, 
an analogue of Theorem
 \ref{thm:Borel-A2} is true
for $\Aut_0(\AA^2)$
after replacing $\JONQ^+(\AA^2)$ by 
$\Jonq^+(\AA^2)$
 and replacing in (d) 
 the optimal upper bound 5 by
the optimal upper 
bound 4.  
\bcor\la{cor:connect-solv-1} 
Let $H\subset \Aut_0(\AA^2)$ be
a solvable connected  subgroup. Then $H$
is elliptic and conjugate in $\Aut_0(\AA^2)$ 
to a subgroup of $\Jonq^+(\AA^2)$. 
Thus, $\Jonq^+(\AA^2)$ is a Borel subgroup
of $\Aut_0(\AA^2)$. 
Any Borel subgroup of $\Aut_0(\AA^2)$
 is maximal among 
the solvable subgroups of $\Aut_0(\AA^2)$
and among the closed proper 
subgroups of $\Aut_0(\AA^2)$,
and coincides with its normalizer 
in $\Aut_0(\AA^2)$.
\ecor
\bproof By Corollary 
\ref{cor:connect-solv}, 
$H$ cannot 
be parabolic, since the latter is 
a torsion subgroup, 
see Theorem \ref{thm:Lamy0}(a). 
Being connected, $H$ cannot be 
elementary loxodromic, 
see Theorem \ref{thm:Lamy0}(b). 
Being solvable, $H$ cannot be 
general loxodromic, 
see Theorem \ref{thm:Lamy0}(d). 
Therefore, $H$ is elliptic. If $H$  
is conjugate to a subgroup of 
 $\GL(2, \kk)$, then it is conjugate to 
a subgroup of $B^+\subset\Jonq^+(\AA^2)$
by the Lie-Kolchin Theorem. 
The  properties of maximality 
and self-normalization
follow from Theorem \ref{thm:BEE}(e) 
due to the decompositions
\[\Aut(\AA^2)=\Trans(\AA^2)
\rtimes{\Aut}_0(\AA^2)
\quad\text{and}\quad {\JONQ}^+(\AA^2)
=\Trans(\AA^2)\rtimes
 {\Jonq}^+(\AA^2),\]
where 
${\Trans}(\AA^2)\subset\Aff(\AA^2)$ 
is the subgroup of translations 
(cf. \cite[Proposition~5(e)]{BEE16}).
\eproof
In the following Theorem 
\ref{prop:max}, we give
examples of 
solvable disconnected algebraic 
subgroups  of $\Aut_0(\AA^2)$ 
and $\Aut(\AA^2)$
which are maximal among 
the solvable subgroups of 
$\Aut_0(\AA^2)$ and 
$\Aut(\AA^2)$, respectively.
Such a subgroup is not contained in 
any Borel subgroup.  
These examples show that  
the connectedness assumption 
in Theorem \ref{thm:BEE}(c) is essential. 
Moreover, we classify all subgroups 
with these properties
up to conjugation. 

Let $Z\subset\TT$ stand for the 
subgroup of scalar matrices and
$2O\subset \SL(2,\kk)$ for
the binary octohedral subgroup
of order 48. Note that 
the subgroup of translations
$\Trans(\AA^2)$ is 
maximal among the abelian subgroups
of $\Aut(\AA^2)$, see  
\cite[Lemma~7.45]{Lam25}.
\bthm\la{prop:max} $\,$
\bnum
\item[$(a)$] The subgroups 
$\langle \TT, \tau\rangle$ 
and $\langle Z, 2O\rangle$ of $\GL(2,\kk)$
are  maximal solvable in $\Aut_0(\AA^2)$. 
Any maximal solvable algebraic
subgroup of $\Aut_0(\AA^2)$
is conjugate to one of them. 

\item[$(b)$] 
The subgroups 
$\langle \TT, \tau, \Trans(\AA^2)\rangle$ 
and
$\langle Z, 2O,  \Trans(\AA^2)\rangle$
of $\Aff(\AA^2)$
are maximal solvable in $\Aut(\AA^2)$.
Any maximal solvable algebraic
subgroup of $\Aut(\AA^2)$ is conjugate to 
one of them.
\enum
\ethm
\bproof
(a) Let us first show  that 
$\langle \TT, \tau\rangle$ and 
$\langle Z, 2O\rangle$
are maximal among the 
solvable subgroups of 
$\Aut_0(\AA^2)$.
Let $H\subset \Aut_0(\AA^2)$ be
a solvable subgroup 
that contains 
$\langle \TT, \tau\rangle$, resp. 
$\langle Z, 2O\rangle$. 
We can suppose that $H$ is closed.
Since $\TT\simeq 
\mathbb{G}^2_{\mathrm m}$ 
and $Z\simeq \Gm$ 
are not torsion subgroups, $H$ 
cannot be parabolic, 
see Theorem \ref{thm:Lamy0}(a).
It also cannot be elementary loxodromic, 
because otherwise $H^0=1$, 
see Theorem \ref{thm:Lamy0}(b), 
which is contradictory. 
Furthermore,
$H$ cannot be focal 
loxodromic and, 
being solvable, it is not 
 general loxodromic, see 
Theorem \ref{thm:Lamy0}(c) and~(d).
Thus, $H$ is elliptic. 
It is therefore conjugate 
to a subgroup, say 
$S$ of either $\GL(2,\kk)$ or 
$\Jonq^+(\AA^2)$. 

Let $H$ contains $\langle \TT, \tau\rangle$.
If $S=H^\phi\subset\GL(2,\kk)$, 
then we can choose
 $\phi\in\GL(2,\kk)$ such that
$\TT^\phi=\TT\subset S^0\subset B^+$,
see Lemma \ref{prop-quasitorus}. 

Let $s\in S\setminus\TT$ 
be the image of $\tau\in H\setminus\TT$.
Then  $s$ normalizes $\TT$, so
$\langle \TT,s\rangle=
\langle \TT,\tau\rangle\subset S$.
By Proposition \ref{prop:GL2}(c) we have
$S=\langle \TT,\tau\rangle=H$, as desired.

Suppose that 
$S\subset \Jonq^+(\AA^2)$. 
The same argument as before gives
$\tau\in S\subset  \Jonq^+(\AA^2)$,
which is a contradiction because
$\tau\notin \Jonq^+(\AA^2)$.

Finally, suppose that $H$ 
contains $Z\times 2O$.
Then the derived length of
$H$ and $S$ is
$\ge 4$.
Thus, $S$ cannot be contained 
in the metabelian group 
$ \Jonq^+(\AA^2)$, and so
$S\subset\GL(2,\kk)$. 
By Proposition  \ref{prop:GL2}(c) 
we have $S=Z\times 2O$. 
This proves that $\langle \TT, \tau\rangle$ 
and $\langle Z, 2O\rangle$
are maximal 
solvable subgroups of 
$\Aut_0(\AA^2)$.

Now, let $R_0$ be a
maximal solvable
algebraic
subgroup of $\Aut_0(\AA^2)$.  
It is  disconnected, 
since otherwise it is properly 
contained in a Borel subgroup, 
and so cannot be maximal solvable.
Being algebraic, $R_0$ 
is elliptic, see Theorem \ref{thm:Lamy0}. 
Being disconnected, $R_0$ 
is not a Borel subgroup of 
$\Aut_0(\AA^2)$. 
Being maximal solvable, $R_0$ 
is not conjugate to 
a subgroup of $\Jonq^+(\AA^2)$.
It is therefore conjugate to a subgroup 
of $\GL(2,\kk)$. Being maximal solvable,
$R_0$ is conjugate to one of the subgroups
$\langle \TT, \tau\rangle$ and
$\langle Z, 2O\rangle$, see Proposition
\ref{prop:GL2}. 

(b) Let $R\subset\Aut(\AA^2)$ 
be a solvable subgroup that contains 
$\langle \TT, \tau, {\rm Trans}(\AA^2)\rangle$.
Then $H:=R\cap  \Aut_0(\AA^2)$ 
is a solvable subgroup of $\Aut_0(\AA^2)$
containing $\langle \TT, \tau\rangle$.
By (a), $H=\langle \TT, \tau\rangle$.
On the other hand,  $ t_{-v}\circ r\in H$ 
for every $r\in R$,
where $v=r(0)$ and $t_{-v}\in {\rm Trans}(\AA^2)$ 
is the translation on $-v$. Thus, 
$R= {\rm Trans}(\AA^2)\rtimes H=
\langle \TT, \tau, {\rm Trans}(\AA^2)\rangle$.
Therefore,  
$\langle \TT, \tau, {\rm Trans}(\AA^2)\rangle$
is a maximal solvable subgroup of $\Aut(\AA^2)$.
The same argument shows that 
$\langle Z, 2O,  \Trans(\AA^2)\rangle$
is also maximal solvable in $\Aut(\AA^2)$.
This proves the first assertion of (b). 
The proof of the second repeats that of (a). 
\end{proof}
\section{Affine toric surfaces and
Borel subgroups}\la{sec:toric}
\subsection{Generalities}
Let $d$ and $e$ be coprime integers 
with $0<e<d$, and let
$\zeta\in\mathbb{G}_{\mathrm m}$ 
be a primitive root of unity 
of degree $d$. Consider
an affine toric surface 
$ X_{d,\,e} \, = \, \AA^2/\, G_{d,\, e}$,
where $G_{d,\, e}=\langle g\rangle\subset\TT$ 
is the cyclic subgroup of order $d$ generated
by $ g=\begin{pmatrix}
\zeta^e & 0 \\
0 & \zeta
\end{pmatrix}\in\TT$.  
\brems\la{rem:conjug} $\,$ 

1. The surface $X_{d,\,e}=
\Spec (k[x,y]^{G_{d,\, e}})$ is 
a toric affine surface 
with the acting torus 
$\TT/G_{d,\, e}\simeq\mathbb{G}_{\mathrm m}^2$, 
cf.~\cite[Section~2.6]{Ful93}. 
It is normal and has
a  unique singular point 
$\pi(0)$, 
where $\pi \colon \AA^2 \to X_{d,\,e}$ is the
quotient morphism and  $0$ 
stands for the origin of $\AA^2$. 

2.
The surface $X_{d,1}$ is 
the Veronese affine cone 
over the rational normal curve 
$\Gamma_d$ of degree $d>1$ in $\PP^d$. 
The group $\GL(2,\CC)$ 
acts naturally on $X_{d,1}$
via the standard irreducible representation 
of $\GL(2,\CC)$
on the space of binary forms of degree $d$.
This action induces an action of 
${\rm PGL}(2, \CC)$ 
on the pair $(\PP^d,\Gamma_d)$.

3. 
The surfaces
$X_{d,\,e}$ and $X_{d',\,e'}$ are isomorphic 
if and only if $d=d'$
and either $e=e'$ or $ee'=1 \mod d$. 
In the latter case, they
are related via the twist 
$\tau:(x,y)\mapsto (y,x)$ 
on $\AA^2$; indeed, 
we have $\tau G_{d,e}\tau=G_{d,e'}$,
where $ee'\equiv 1 \mod d$.
On the other hand, given $X_{d,\,e}$
we have ${\rm Cl}(X_{d,e})=\ZZ/d\ZZ$, and
the integer $e$ can be reconstructed 
from the minimal 
resolution of singularity of $X_{d,e}$ up to
inversion $e\mapsto e'\equiv e^{-1}\mod d$,
see \cite[P.~47, last exercise, 
and P.~136, item 23]{Ful93}.

Furthermore, $\tau$ normalizes 
$G_{d,e}$ if and only if $e^2\equiv 1\mod d$.
In the latter case, $\tau$ acts on $G_{d,e}$ 
via $g\mapsto g^e$.
\erems
\bnota\la{not:N} We let
\[N_{d,e}={\Norm}_{\GL(2,\kk)}
(G_{d,\,e})\quad\text{and}
\quad\cN_{d,\,e}=
{\Norm}_{\Aut(\AA^2)}(G_{d,\,e}).\]
It is easily seen that
\[
N_{d,e}=\begin{cases} \GL(2,\kk) & \text{ if   
$e=1$,}\\
\Norm_{\GL(2,\kk)}(\TT)=
\langle\TT,\tau\rangle  &\text{ if  
$ e\neq 1$  and
$e^2\equiv 1\mod d$,}\\
\TT&\text{ if  
$e^2\not\equiv 1\mod d$.}\end{cases}
\]
We let $B^\pm\subset\GL(2,\kk)$ denote 
the Borel subgroup of all upper
(lower, respectively) triangular matrices.
Consider the  subgroups 
\be\label{nde} 
N_{d,e}^\pm=N_{d,e}\cap
{\Jonq}^{\pm}(\AA^2)=\begin{cases}
B^\pm &\text{ if $e=1$,}\\
\TT&\text{ otherwise}
\end{cases}
\ee 
and
 \[ \cN_{d,e}^\pm=\cN_{d,e}\cap 
 {\Jonq}^{\pm}(\AA^2)\,.\]
Note that \be\label{tor-int}
\cN_{d,e}^+\cap\cN_{d,e}^-= 
{\Jonq}^{+}(\AA^2) \cap
{\Jonq}^{-}(\AA^2)=\TT\,.\ee
\enota
\blem[{\rm \cite[Lemmas~4.5--4.6]{AZ13}}]
\la{lem-descr-normalis-bis} $\,$ 
Consider the $\ZZ/d\ZZ$-grading
\be\la{eq:grading}
\kk[t]=\bigoplus_{i=0}^{d-1} A_{d,i},
\quad\text{where}\quad A_{d,i}=t^i\kk[t^d].\ee
\bnum\item[$(a)$] 
The group
$\cN_{d,e}^+$ ($\cN_{d,e}^-$, respectively) 
consists of all de
Jonqi\`eres transformations $\varphi^+$ 
as in (\ref{eq-Jonq-plus})
($\varphi^-$ as in (\ref{eq-Jonq-minus}), 
respectively) with $tp(t)\in
A_{d,e}$ ($tp(t)\in A_{d,e'}$, respectively).
\item[$(b)$] The subgroup $\cN_{d,e}^\pm$
coincides with the centralizer 
$\Cent_{\Jonq^\pm(\AA^2)}(G_{d,e})$.
\item[$(c)$] For all $\varphi\in\cN_{d,e}$ 
we have $\varphi(0)=0$.
\item[$(d)$] 
We have $\deg(\varphi)\ge e$ 
if $\varphi\in\cN^+_{d,e}\setminus\TT$ 
and $\deg(\varphi)\ge e'$ 
if $\varphi\in\cN^-_{d,e}\setminus\TT$, 
where $1\le e'< d$ 
and $ee'\equiv 1 \mod d$.
\enum 
\elem
\brems\label{rem:connected} $\,$

1.
By Lemma \ref{lem-descr-normalis-bis}(c) 
the $\cN_{d,e}$-action on $\AA^2$ 
fixes the origin. Respectively, 
any automorphism of the  
affine surface $X_{d,e}$ 
fixes the unique singular point 
$\pi(0)\in X_{d,e}$, where
$\pi\colon\AA^2\to X_{d,e}=\AA^2/G_{d,e}$ 
is the natural surjection. 

2. 
We have $\cN_{d,e}\cap\GL(2,\kk)=N_{d,e}$.
In particular, $\deg\varphi>1$ 
$\forall \varphi\in\cN_{d,e}\setminus N_{d,e}$.
The tangent representation 
of $\cN_{d,e}$ at the origin $0\in\AA^2$ 
is given by $N_{d,e}$. For $e>1$ we have 
$d\varphi(0)\in \langle\TT,\tau\rangle$ 
if $e^2\equiv 1\mod d$ and 
 $d\varphi(0)\in \TT$ if $e^2\not\equiv 1\mod d$,
see Notation \ref{not:N}.

3. 
 Consider the closed ind-subgroup 
$\cA^{d,e}\subset {\Jonq}^+(\AA^2)$
consisting of all $\varphi^+$ 
as in \eqref{eq-Jonq-plus}
with \[tp(t)\in A^{d,e}:=t^d\kk[t^d]\oplus
\bigoplus_{i=1,\ldots,d-1,\,\,i\neq e} A_{d,i},\]
see  \eqref{eq:grading}.
Since $t\kk[t]=A_{d,e}\oplus A^{d,e}$ 
we have a splitting
(see Lemma \ref{lem-descr-normalis-bis}(a))
\be\la{eq:split} {\Jonq}^+(\AA^2)=\cN_{d,e}^+
\ltimes \cA^{d,e}. \ee
\erems
 \bthm[{\rm\cite[Theorem~4.2, Lemma~4.3 
 and Proposition~4.4]{AZ13}
 \footnote{See also \cite[Remark~4.6]{Kov15} 
 for an alternative proof of (b) and (c).}}]
 \la{thm-JvdK-toric-surf} $\,$
 \bnum
 \item[$(a)$] There is an isomorphism
\[\Aut (X_{d,e}) \simeq \cN_{d,e}/G_{d,e}.\]
\item[$(b)$]  For $e^2\not\equiv 1 \mod d$, 
we have
\be\la{eq-amalg-1} 
\cN_{d,e} \simeq \cN_{d,e}^+*_{\TT}
\cN_{d,e}^- \quad\text{and}\quad
\Aut (X_{d,e}) \simeq
\cN_{d,e}^+/G_{d,e}*_{\TT/G_{d,e}} 
\cN_{d,e}^-/G_{d,e}.
\ee
\item[$(c)$] For $e^2\equiv 1 \mod d$, we have 
\be\la{eq-amalg-2} 
\cN_{d,e}  \simeq {\cN_{d,e}^+}*_{N_{d,e}^+}
N_{d,e} \quad\text{and}\quad
\Aut(X_{d,e}) \simeq 
\cN_{d,e}^+/G_{d,e}*_{N_{d,e}^+/G_{d,e}}
N_{d,e}/G_{d,e}.
 \ee 
 \enum
\ethm
\brems \la{rem:2-comp} $\,$

1. Note that the subgroups
$\cN^\pm_{d,e}\subset \Jonq^{\pm}(\AA^2)$ 
are solvable. 

2. It follows from Lemma 
\ref{lem-descr-normalis-bis}(a) and (b)
that the nested ind-subgroups 
$\cN^\pm_{d,e}\subset \Aut(\AA^2)$
are connected and closed. 
The normalizer 
$\cN_{d,e}$ of $G_{d,e}$
is also closed. 
The subgroup
$\langle \cN_{d,e}^+,\cN_{d,e}^-\rangle$ 
of $\cN_{d,e}$
 is connected and closed. 
 It coincides with $\cN_{d,e}$ 
unless $e^2\equiv 1\mod d$. 
In the latter case, 
its index is $2$ in $\cN_{d,e}$ and 
$\cN_{d,e}=\langle \cN_{d,e}^+
,\cN_{d,e}^-, \tau\rangle$, 
see \cite[Lemma~4.10]{AZ13}.
In the case where $e>1$ and 
$e^2\equiv 1\mod d$, we have 
\be\la{eq:decomp}
\cN_{d,e} = {\cN_{d,e}^+}*_{\TT} 
\langle\TT,\tau\rangle=\{g=a_1\tau a_2\tau
\cdots a_{n-1}\tau  a_n\},\,\,\, a_i\in
\cN_{d,e}^+\,\,\,\text{and}\,\,\, a_i\notin\TT\,\,\,
\text{for}\,\,\, 1<i<n.\ee
Since $\cN_{d,e}^+$ is connected, 
in the latter case, $\cN_{d,e}$ 
has two connected components.
The unity component consists of 
all $g$ with $n$ odd, while
the other component consists of 
all $g$ with $n$ even. 
Another proof of this 
statement is given in
\cite[Proposition~3]{Kik24}.
\erems
By Serre's theorem, an analogue 
of Theorem \ref{Serre} 
holds for the groups  $\cN_{d,e}$ and 
$\Aut(X_{d,e})=\cN_{d,e}/G_{d,e}$, 
see \cite[Theorem~4.15]{AZ13}. 
Furthermore, 
the following analogue of Lemma 
\ref{prop-quasitorus} holds, see 
\cite[Corollary~4.16, Theorem~4.17 
and its proof]{AZ13}.
\bprop $\,$
\bnum
\item[$(a)$]
Any reductive algebraic subgroup $G\subset
\Aut(X_{d,e})$ is conjugate to 
a subgroup of $N_{d,e}/G_{d,e}$. 
Furthermore, 
if $G\subset \cN_{d,e}^\pm/G_{d,e}$, then
$G^{\pi(\mu)}\subset N_{d,e}^\pm/G_{d,e}$ 
for some
$\mu\in U^\pm\cap \cN_{d,e}^\pm$, 
where $\pi\colon \cN_{d,e}\to \cN_{d,e}/G_{d,e}$ 
is the natural homomorphism. 
\item[$(b)$] Every unipotent algebraic subgroup 
of $\Aut(X_{d,e})$ is commutative.
\enum
\eprop
\subsection{Borel subgroups of $\Aut(X_{d,e})$}
Recall that
a \emph{Borel subgroup} of an ind-group $G$
 is a solvable connected subgroup 
that is maximal under inclusion among the solvable  
connected subgroups of $G$, 
see Definition \ref{def:Borel}. 
\bprop \la{prop:Borel-Nde} The subgroup
$\cN^\pm_{d,e}$ (resp., $\cN^\pm_{d,e}/G_{d,e}$)
is a Borel subgroup of $\cN_{d,e}$ 
(resp., of $\Aut(X_{d,e})=\cN_{d,e}/G_{d,e}$). It is
maximal among the solvable subgroups 
and among the  closed proper subgroups 
of $\cN_{d,e}$ 
(resp., of $\Aut(X_{d,e})$). In particular, 
$\cN^\pm_{d,e}$ 
(resp., $\cN^\pm_{d,e}/G_{d,e}$) 
coincides with its normalizer.
\eprop
\bproof
According to Remarks
\ref{rem:2-comp},
the $\cN^\pm_{d,e}$ 
are solvable connected  subgroups
 of $\cN_{d,e}$. 
 The same is true for the quotient groups 
 $\cN^\pm_{d,e}/G_{d,e}$. 
 Applying Proposition \ref{thm:maximal}
to the amalgamated products in
\eqref{eq-amalg-1} and \eqref{eq-amalg-2}
we obtain the maximality among 
the solvable subgroups. 
Due to the splitting \eqref{eq:split},
the maximality of $\cN^+_{d,e}$
among 
the closed proper subgroups follows from
that for $\Jonq^+(\AA^2)$, 
see Corollary \ref{cor:connect-solv-1}.
The same argument works for $\cN^-_{d,e}$.
This maximality survives 
when passing to the quotient, 
and so it also holds for 
$\cN^\pm_{d,e}/G_{d,e}$.
\eproof
\bprop\la{prop:maxim}
Any solvable subgroup 
$G\subset\cN_{d,e}$ 
(resp. $G\subset\cN_{d,e}/G_{d,e}$) 
is contained in
 a maximal solvable subgroup of $\cN_{d,e}$ 
 (resp. of $\cN_{d,e}/G_{d,e}$).
Any connected 
solvable subgroup $G\subset\cN_{d,e}$ 
(resp. $G\subset\cN_{d,e}/G_{d,e}$) 
is contained in
a Borel subgroup of $\cN_{d,e}$ 
(resp. of $\cN_{d,e}/G_{d,e}$).
\eprop
\bproof The conclusions follow immediately from 
Proposition \ref{prop:FP}.
Indeed, the assumption of the latter 
proposition is fulfilled in both cases,
since the derived length of solvable subgroups of
$\cN_{d,e}$ and of $\cN_{d,e}/G_{d,e}$ 
does not exceed $5$, see
Theorem \ref{thm:BEE}(d). 
\eproof
The next theorem is one of the main 
ingredients in the proof of 
Theorem \ref{main} from 
the Introduction. 
\bthm\label{thm:Borel} $\,$
\bnum
\item[$(a)$] 
Every Borel subgroup of $\cN_{d,e}$
is conjugate in $\cN_{d,e}$  either to 
$\cN^+_{d,e}$  or to $\cN^-_{d,e}$. 

\item[$(b)$] Every Borel subgroup of
$\Aut(X_{d,e})
= \cN_{d,e}/G_{d,e}$ is
conjugate in 
$ \cN_{d,e}/G_{d,e}$  either to
$\cN^+_{d,e}/G_{d,e}$ or to
$\cN^-_{d,e}/G_{d,e}$.
\enum
\ethm
\brem \la{rem:tau}
If $e^2\equiv 1\mod d$, then $\tau\in\cN_{d,e}$
and $\cN^-_{d,e}=\tau\cN^+_{d,e}\tau$ 
is conjugate to $\cN^+_{d,e}$. Otherwise
$\tau\notin\cN_{d,e}$ and the factors
$\cN^+_{d,e}$ and $\cN^-_{d,e}$ 
in \eqref{eq-amalg-1}
are not conjugate in $\cN_{d,e}$.
Indeed, the conjugacy of the factors 
$A$ and $B$
in a proper amalgam $G=A*_{A\cap B} B$ 
would contradict the uniqueness 
of the reduced word 
presentation.
\erem
The proof of Theorem \ref{thm:Borel} 
is preceded by the following 
Lemma \ref{lem:finite-geod}. 
Its proof is similar to the proof 
of \cite[Lemma~7.50(2)]{Lam25}. 
For the reader's convenience, we provide 
an argument. 

\smallskip

We let $T_{d,e}$ 
be the Bass--Serre tree 
of the amalgam 
$\cN_{d,e}
=A*_C B$
as in  \eqref{eq-amalg-1} 
resp. \eqref{eq-amalg-2}.

\blem \la{lem:finite-geod}  
Assume that $e^2\not\equiv 1\mod d$.
Let  $\gamma$ be a geodesic segment  
of $T_{d,e}$
of length $3$ and
$\Stab(\gamma)\subset\cN_{d,e}$
the subgroup consisting
of all $f\in\cN_{d,e}$
fixing the vertices of $\gamma$.
Then the following hold.
\bnum
\item[$(a)$] Every
$f\in\Stab(\gamma)$ is 
conjugate in $\cN_{d,e}$ to 
$f'\colon (x,y)\mapsto (\alpha x, \beta y)$,
where $\alpha$ and $\beta$ 
are  primitive roots of unity
 of the same degree;
\item[$(b)$] $\Stab(\gamma)$
is conjugate 
in $\cN_{d,e}$ to 
$G_{d',e'}$ 
for some pair $(d',e')$
of coprime integers;
\item[$(c)$] 
there is a loxodromic $g\in\cN_{d,e}$ 
that centralizes $\Stab(\gamma)$;
\item[$(d)$]  for every  
$f\in \Stab(\gamma)$ 
the subtree of fixed points $T_{d,e}^f$ 
contains 
the geodesic $\Ax(g)$,  and therefore
is unbounded.
\enum
\elem
\bproof (a)
Note that $\cN_{d,e}$ acts transitively 
on the set of edges
of $T_{d,e}$. Hence, up to conjugation 
we may 
assume that $\gamma$ includes the edge 
$[1\cdot A-1\cdot B]$ 
as the middle edge, where the dash
represents as usual the edge 
$1\cdot (A\cap B)=1\cdot \TT$. 
Thus, for some 
$a\in A\setminus B$ and $b\in B\setminus A$,
the geodesic segment $\gamma$ 
takes the form 
\be\la{eq:gamma}
\gamma=
[b\cdot A-1\cdot B-1\cdot A-a\cdot B]\ee
where the dashes represent
the edges $b\cdot \TT, 1\cdot \TT$, 
and $a\cdot \TT$, 
respectively.

Recall that
$\cN^\pm_{d,e}=
\Cent_{\Jonq^\pm}(G_{d,e})$, 
see Lemma \ref{lem-descr-normalis-bis}(b).
Knowing that   $A\cap B=\TT$,  we can
choose 
$a\in A\setminus B$ 
and $b\in B\setminus A$ of the form 
\be\la{eq:lin-parts}
a\colon (x,y)\mapsto (x+y^{e}p(y^d),y)
\quad\text{and}\quad
b\colon (x,y)\mapsto (x,y+x^{e'}q(x^d)),
\quad\text{where}\quad 
p,q\in\kk[t]\setminus\{0\},\ee
see Corollary \ref{cor:commut}. 
Since $f$ 
fixes the vertices $1\cdot A$ and $1\cdot B$ 
of $T_{d,e}$, we have
$f\in A\cap B=\TT$, and therefore, 
$f\colon (x,y)\mapsto (\alpha x, \beta y)$ 
for some 
$\alpha, \beta\in\mathbb{G}_{\mathrm m}$. 
Since $f$ 
fixes the vertices $bA$ and $aB$ we have
\[b^{-1}fb=a_1\in A\cap B=
\TT\quad\text{and}
\quad a^{-1}fa=b_1\in B\cap A=\TT.\]
Furthermore, $f$ and $a_1$ 
(resp. $f$ and $b_1$) are conjugate 
in $\GL(2,\kk)$ 
via $da|_{T_0\AA^2}=1$ 
(resp.  $db|_{T_0\AA^2}=1$), 
see  \eqref{eq:lin-parts}.
It follows that $a_1=b_1=f$, 
and therefore $f$ commutes with
$a$ and $b$. 

These commutator relations
lead to the equalities
\be\la{eq:compar} \alpha p(y^d)
=\beta^ep(\beta^dy^d)  
\quad\text{and}\quad \beta
q(x^d)=\alpha^{e'}q(\alpha^dx^d),
\quad\text{where}\quad e,e'>1.\ee
Comparing the coefficients of the 
monomials of lowest degree 
in $p$ and $q$ on either side 
of \eqref{eq:compar}, 
we obtain the relations
\[\alpha=\beta^{e+kd}
\quad\text{and}\quad
 \beta=\alpha^{e'+ld},
 \quad\text{and so}\quad 
 \alpha^{md}=\beta^{md}=1,\]
where 
 $m=(e+kd)(e'+ld)-1\equiv 
 0 \mod d$.
 It follows that $\alpha$ and $\beta$ 
 are primitive roots of unity 
 of the same degree
 that divides $md$. 
 This proves (a).
  
(b) Note that the integers $k,l$ 
and $m$ 
depends only on $a$ and $b$, 
and do not depend on 
$f\in\Stab(\gamma)$.
Since $m$ and $e+kd$
are comprime and $d$ divides $m$,
we have $\gcd(e+kd, md)=1$. 
Let $\tilde d=md$ and 
$\tilde e\in[1,\tilde d-1]$ 
be the residue 
of $e+kd$ modulo $\tilde d$.
Then
every $f=(\alpha x,\beta y)\in
\Stab(\gamma)$
belongs to $G_{\tilde d,\tilde e}$,
see the proof of (a). 
Therefore, $\Stab(\gamma)\subset 
G_{\tilde d,\tilde e}$ is a cyclic group
that coincides with some 
$G_{d',e'}$, where $d'|\tilde d$.
Now (b) follows.

(c)-(d) Since $f\in\Stab(\gamma)$
commutes with $a$ and $b$, it 
also commutes with the loxodromic 
element $g=ab\in\cN_{d,e}$.
Hence $T_{d,e}^f\supset \Ax(g)$ 
is unbounded, 
see Proposition \ref{prop:elements}(b). 
\eproof
\brem Let $\gamma$ be a segment in
\eqref{eq:gamma}. 
Then the loxodromic elements 
$ab, ba\in \cN_{d,e}$ 
centralize $\Stab(\gamma)$, 
and for their geodesic axes we have 
\[\Ax(ab)=] \ldots aba\cdot B-ab\cdot A-
a\cdot B-1\cdot A-1\cdot B-b^{-1}\cdot A-
b^{-1}a^{-1} B-b^{-1}a^{-1}b^{-1}\cdot A
 \ldots[\]
resp.
\[\Ax(ba)=] \ldots a^{-1}b^{-1}a^{-1}\cdot B-
a^{-1}b^{-1}\cdot A-a^{-1}\cdot B-1\cdot A-
1\cdot B-b\cdot A-ba\cdot B-bab\cdot A 
\ldots[\]
We see that
\[\Ax(ab)\cap\gamma=
[a\cdot B-1\cdot A-1\cdot B],\,\,\, 
\Ax(ba)\cap\gamma=
[1\cdot A-1\cdot B-b\cdot A],
\,\,\,
\Ax(ab)\cap\Ax(ba)
=[1\cdot A-1\cdot B].\]
\erem
\bcor\la{cor:diam} The groups
$\cN_{d,e}$ and 
$\Aut(X_{d,e})=\cN_{d,e}/G_{d,e}$ 
are torsionally unbounded.
\ecor
\bproof 
It suffices to prove the assertion for 
$\cN_{d,e}$. Indeed, it is also valid for
$\cN_{d,e}/G_{d,e}$ thanks to 
Lemma \ref{lem:quotient}.

Let us first assume that
$e^2\equiv 1\mod d$. 
In this case, 
$\cN_{d,e}  \simeq 
{\cN_{d,e}^+}*_{N_{d,e}^+}
N_{d,e} $ inherits 
the amalgam structure 
from that of $\Aut_0(\AA^2)$, 
see \eqref{eq-amalg-2} and 
Lemma \ref{lem:amalg-A0}. 
Therefore, there is 
an embedding of trees
$T_{d,e}:=T(\cN_{d,e})\subset T_0
=T(\Aut_0(\AA^2))$, 
see Remark \ref{rem:tree-embed}.
Thus, $\cN_{d,e}$ acts on $T_0$ 
leaving 
the subtree $T_{d,e}$ invariant. 
So, the assertion follows from 
Proposition \ref{prop:unbound}. 

Now let $e^2\not\equiv 1\mod d$. 
Then
for every elliptic non-torsion element 
$f\in\cN_{d,e}$  
the diameter of the fixed subtree 
$T^f_{d,e}$ 
does not exceed $2$. 
Indeed, if ${\rm diam}(T^f_{d,e})\ge 3$, 
then 
$T^f_{d,e}$ contains a geodesic segment 
of length  $3$.
By Lemma \ref{lem:finite-geod}, 
in this case, $f$ 
is a torsion element. So, once again
$\cN_{d,e}$ is torsionally unbounded.
\eproof
\bproof[Proof of Theorem \ref{thm:Borel}]
We provide a proof of (a); 
the proof of (b) is identical. 
By Corollary~\ref{cor:diam}, $\cN_{d,e}$  
acting on 
the Bass-Serre tree $T_{d,e}$ 
is torsionally unbounded. 
Therefore, by Corollary~\ref{cor:clef}
every Borel subgroup $\cB$ 
of $\cN_{d,e}$ is elliptic
and is conjugate to a Borel subgroup of
one of the factors $A$ and $B$ 
of the amalgam 
$\cN_{d,e}=A_{A\cap B} B$
 in \eqref{eq-amalg-1} 
or \eqref{eq-amalg-2}.
Recall that
$\cN_{d,e}^\pm$ are Borel subgroups of 
$\cN_{d,e}$, see Proposition  
\ref{prop:Borel-Nde}. 

If $e^2\not\equiv 1\mod d$, 
then $A=\cN_{d,e}^+$ and 
$B=\cN_{d,e}^-$, see \eqref{eq-amalg-1}.
Therefore, $\cB$ is conjugate 
to one of these
subgroups.

Let $e=1$.  Then
$\cN_{d,e}=\cN^+_{d,e}*_{N^+} \GL(2,\kk)$.
If $\cB$ is conjugate 
to a Borel subgroup of 
$\GL(2,\kk)$, 
then it is also conjugate 
to the triangular subgroup
$B^+\subset \cN_{d,e}^+$ by 
the Lie-Kolchin theorem. 
Since $\cB$ is a Borel subgroup of 
$\cN_{d,e}$
and the inclusion above is strict, 
this leads to a contradiction.
We conclude that $\cB$ is conjugate to 
$\cN_{d,e}^+$.

Let $e>1$ and $e^2\equiv 1\mod d$.
Then 
$\cN_{d,e}=\cN^+_{d,e}*_{\TT} N_{d,e}$,
where $N_{d,e}=\langle\TT,\tau\rangle$.
Thus, $\cB$ is conjugate either 
to $\cN_{d,e}^+$ or 
to the Borel subgroup $\TT$
 of $N_{d,e}$.
The latter case
is impossible, because the inclusion  
$\TT\subset  \cN_{d,e}^+$ is strict.
\eproof
The following theorem is an analogue 
of \cite[Theorem~4]{BEE16}.
\bthm\la{thm:no-fin-ind} 
A maximal solvable non-abelian subgroup 
$H$ of 
$\cN_{d,e}$ (resp., of $\Aut(X_{d,e})$) 
is a Borel subgroup 
if and only if it contains no proper 
 subgroup of finite index.
\ethm
The proof is preceded by Lemmas 
\ref{lem:only-if} and \ref{lem:ell}. 
\blem\la{lem:only-if} A Borel subgroup $\cB$ of 
$\cN_{d,e}$ (resp., of $\Aut(X_{d,e})$) 
is a non-abelian 
maximal solvable subgroup
that contains no proper 
subgroup of finite index.
\elem
\bproof 
By Theorem \ref{thm:Borel}, $\cB$ 
is conjugate to one of  the
$\cN_{d,e}^\pm$
(resp., to one of the $\cN_{d,e}^\pm/G_{d,e}$). 
Therefore, $\cB$ is non-abelian, 
connected and nested. 
By Proposition \ref{prop:Borel-Nde}, 
$\cB$ is maximal solvable. 
By Lemma \ref{lem:divisible-1}, $\cB$ 
has no proper subgroup
of finite index.
\eproof
\blem\la{lem:ell} Let $H\subset \cN_{d,e}$ 
(resp., 
$H\subset \cN_{d,e}/G_{d,e}$) 
be a maximal 
solvable non-abelian,  elliptic subgroup 
containing no proper 
subgroup of finite index. 
Then $H$ is a Borel subgroup 
of $\cN_{d,e}$ (resp., of 
$\cN_{d,e}/G_{d,e}$).
\elem
\bproof We give a proof 
for the case where
$H\subset \cN_{d,e}$;
 the same proof  works 
 for the case
where 
$H\subset \cN_{d,e}/G_{d,e}$.
If $e^2\not\equiv 1\mod d$,
then, 
by \eqref{eq-amalg-1} 
(see Theorem 
\ref{thm-JvdK-toric-surf}),
$H^\phi\subset\cN_{d,e}^\pm$ 
for some $\phi\in\cN_{d,e}$ 
because $H$ is elliptic.
In fact, this inclusion 
is an equality 
due to the maximality 
of~$H$. 

If $e=1$, then  by
\eqref{eq-amalg-2}, $H$ 
is conjugate either to a subgroup 
of $\cN_{d,e}^+$,
or to a subgroup of $\GL(2,\kk)$. 
In the first case, 
$H=\cN_{d,e}^+$ is a 
Borel subgroup 
of $\cN_{d,e}$ 
due to the maximality 
of $H$. In the second case, $H$ is 
a solvable algebraic subgroup 
of $\GL(2,\kk)$. 
Since $H$ has no proper
subgroup of finite index, 
$H$ is  connected by 
Lemma \ref{lem:alg-grp}.
Due to the maximality of $H$, 
it is a Borel subgroup 
of $\GL(2,\kk)$  conjugate to 
$B^+\subset \cN_{d,e}^+$.
This contradicts the 
assumption that $H$ 
is maximal among the 
solvable subgroups 
of $\cN_{d,e}$. 

If $e>1$ and $e^2\equiv 1\mod d$,
then by
\eqref{eq-amalg-2},
$H$ is conjugate either to 
a subgroup of $\cN_{d,e}^+$ 
or to 
a subgroup of $\langle\TT,\tau\rangle$.
In the first case, $H=\cN_{d,e}^+$ is a Borel 
subgroup  of $\cN_{d,e}$ due 
to the maximality of $H$. 
In the second case, $H$ is conjugate to 
$\TT$, 
because it does not have 
a proper subgroup 
of finite 
index. Since $\TT\subset\cN_{d,e}^+$,
the latter contradicts 
the maximality of $H$. 
\eproof
\bproof[Proof of Theorem 
\ref{thm:no-fin-ind}] The ``only-if'' part 
follows from 
Lemma \ref{lem:only-if}. 
To prove the ``if'' part, suppose that $H$
contains no proper subgroup of finite index. 
Since $H$ is not abelian, 
it cannot be parabolic
(see Theorem 4.13(a)). 
Nor can it be elementary loxodromic, 
since the latter subgroups contain 
proper finite index subgroups
(see Theorem 4.13(b); 
cf. \cite[Proposition~4]{BEE16}). 
According to Theorem 4.13(c), 
$H$ is not focal loxodromic. 
 It is also not generally loxodromic, 
since it is solvable
(see Theorem 4.13(d)). 
Thus, $H$ is elliptic. 
The proof is now complete by applying 
Lemma \ref{lem:ell}. 
\eproof
\bcor Any automorphism 
of the group $\Aut(X_{d,e})$ 
considered as an abstract group sends 
a Borel subgroup of $\Aut(X_{d,e})$ 
to a Borel subgroup. 
\ecor
\brems

1. We know that any automorphism 
of the group $\Aut(\AA^2_{\CC})$ 
considered as an abstract group is inner 
up to an automorphism of $\CC$ over $\QQ$, 
see \cite{Des06}. On the other hand, 
for $e^2\not\equiv 1\mod d$, there exist 
automorphisms of the abstract group 
$\Aut(X_{d,e})$ 
that are not inner up to an automorphism 
of the base  field $\kk$, 
see \cite[Remark~5.13]{LRU23}. 
The automorphisms constructed 
in [ibid] preserve the Borel subgroups 
$\cN_{d,e}^\pm/G_{d,e}$, 
factors of the corresponding 
amalgam in \eqref{eq-amalg-1}.
These automorphisms 
are identical on $\cN_{d,e}^+/G_{d,e}$
and non-trivial on $\cN_{d,e}^-/G_{d,e}$. 

2. Conjugation with involution 
$\tau$ sends 
$G_{d,e}$ to $G_{d,e'}$, 
$\Aut(X_{d,e})$ to $\Aut(X_{d,e'})$,
$\cB^+_{d,e}$  to $\cB^-_{d,e'}$ 
and $\cB^-_{d,e}$ to $\cB^+_{d,e'}$, 
see Remarks \ref{rem:conjug}.3
and \ref{rem:tau}. 
In the case where 
$e^2\equiv 1\mod d$ 
this defines an inner
automorphism of $\Aut(X_{d,e})$
 that interchanges 
 the Borel subgroups
 $\cB^+_{d,e}$ and $\cB^-_{d,e}$.
\erems
\section{Appendix A: Tree amalgams}
\la{sec:tree-amalg}
In this appendix we extend 
Proposition \ref{thm:maximal} 
to amalgams 
of a tree of groups. We also present
examples of bearable groups of automorphisms; 
such a group is an amalgam of a tree of ind-groups. 
\subsection{Amalgamated tree products}
Recall the following definitions.
\bdefi[{\rm cf.~\cite[Introduction and 
Theorem~1]{KS70} and \cite[Definition~8]{Ser03}}]
\la{def:tree-prod}
Let $\Gamma$ be a tree and $(\Gamma, \cG)$
a system  of vertex groups 
$G_v$, $v \in \vert(\Gamma)$, and edge groups 
$G_e$, $e \in {\rm edge}(\Gamma)$, such that
for $e=[u,v]$ there are  proper embeddings 
$G_e\hookrightarrow G_u$
and $G_e\hookrightarrow G_v$. Then there is a 
group $G$
that contains the
$G_u$ and $G_e$ as subgroups and 
verifies the following: 
\begin{itemize}
\item[$(i)$] $G=
\langle G_u\,|\,u\in\vert(\Gamma)\rangle$;
\item [$(ii)$]  for every edge 
$e=[u,v] \in {\rm edge}(\Gamma)$
we have $G_e=G_u\cap G_v$ and
$ \langle G_u,G_v\rangle\simeq 
G_u*_{G_e} G_v$,
where the inclusions
$G_e\subset G_u$ and $ G_e\subset G_v$
correspond to the given embeddings 
$G_e\hookrightarrow G_u$
and $G_e\hookrightarrow G_v$;
\item[$(iii)$]  every group that satisfies (i) 
and (ii)  is a quotient of $G$. 
\end{itemize}
This group $G$ is defined uniquely 
up to isomorphism. It is called
an \emph{(amalgamated) tree product of} 
$(\Gamma, \cG)$, or an 
\emph{amalgam of a tree of groups} 
$(\Gamma, \cG)$,
and is denoted 
by $G=\varinjlim(\Gamma, \cG)$, 
see \cite{Ser03}.
\edefi
\brems\la{rem:subtree} $\,$

1. 
For any subtree $\Gamma'$ of 
$\Gamma$, 
the subgroup
$ G'=\langle G_u\,|\,u\in\vert(\Gamma')\rangle$ is 
a tree product 
$G'=\varinjlim(\Gamma', \cG|_{\Gamma'})$, 
see \cite[Theorem~1, (1) and (2)]{KS70}. 

2. A tree product $G=\varinjlim(\Gamma, \cG)$
is unique up to isomorphism. A presentation of $G$ 
can be obtained as a union of the presentations 
of $G_u*_{G_e} G_v$ 
for all edges 
$e=[u,v]\in{\rm edge}(\Gamma)$, see 
\cite[Introduction]{KS70}.
\erems
\bdefi[\emph{contraction of a tree}, 
{\rm see \cite[p.~231, Definition]{KS70}}]
Let a tree $\Gamma$ be partitioned into  
subtrees $\Gamma_\alpha$ 
which are "disjoint" in the sense that
each vertex of $\Gamma$ belongs to exactly
one of the $\Gamma_\alpha$. 
A \emph{contraction of $\Gamma$
according to this partition} is the tree 
$\gamma$ with vertices
$v_\alpha=\Gamma_\alpha$
and with edges that are the edges of $\Gamma$ 
which join vertices of $\Gamma_\alpha$ 
and $\Gamma_\beta$.
\edefi
\bprop[{\rm \cite[Theorem~1(3)]{KS70}}]\la{prop:KS}
Consider a tree product $G=\varinjlim(\Gamma, \cG)$.
Assume that $\Gamma$ is partitioned into disjoint
subtrees $\Gamma_\alpha$, and let $\gamma$ be
the contraction of $\Gamma$ according to this partition.
 Let also  $G_\alpha$ be the tree product
$G_\alpha=\langle G_v\,|\,v\in\vert(\Gamma_\alpha)\rangle$,
see Remark \ref{rem:subtree}.
Then $G$ is the tree product 
\[G=\varinjlim(\gamma, \cG_\gamma)=\langle G_\alpha\,|\,
v_\alpha\in\vert(\gamma)\rangle,\] where $\cG_\gamma$
has the vertex groups $G_\alpha$ for each
$v_\alpha\in\vert(\gamma)$
and the edge groups $G_\alpha\cap G_\beta$ 
for each $[v_\alpha,v_\beta]\in {\rm edge}(\gamma)$.
\eprop
The following problem arises.

\medskip

\noindent {\bf Problem.}  \emph{Given a tree product 
$G=\varinjlim(\Gamma, \cG)$,
describe the maximal solvable subgroups~of~$G$. }

\medskip

Using Propositions \ref{thm:maximal} 
and \ref{prop:KS}
we deduce the following partial result. 
\bprop Let $G=\varinjlim(\Gamma, \cG)$
be a tree product. Assume that for a vertex 
$v\in\vert(\Gamma)$, the vertex group
$G_v$ is solvable and  contains a subgroup $H_v$
of index $|G_v:H_v|\ge 3$
such that $H_v=G_e$ for any edge $e=[u,v]$ 
of $\Gamma$. 
Then the subgroup $G_v$ is maximal among 
the solvable subgroups of $G$.
\eprop
\bproof Let $\Gamma\ominus\{v\}$ 
be the forest obtained from $\Gamma$
after deleting the vertex $v$ and all 
its incident edges.  
Consider the partition of $\Gamma$ 
into the subtrees
$\Gamma_v=\{v\}$ and $\Gamma_u$, 
where $[u,v]\in {\rm edge}(\Gamma)$
and $\Gamma_u$ is the connected component 
of $\Gamma\ominus\{v\}$
that contains $u$, 
that is, the branch of $\Gamma$ at $v$
containing $u$. 

Consider the graph $\gamma$ 
obtained by contracting $\Gamma$ 
according to this partition. 
Let 
\[A=G_v\quad\text{and}\quad
B={\prod}^*_{[u,v] \in {\rm edge}(\Gamma)} 
B_u,\quad\text{where }\quad
B_u=\varinjlim(\Gamma_u, \cG|_{\Gamma_u})\]
and $\prod^*$ stands for the free product. 
 Note that
the forest $\Gamma\ominus\{v\}$ 
is a disjoint union of 
the branches $\Gamma_u$, and 
$\gamma\ominus\{v\}$
is a disjoint union of one-vertex graphs 
indexed by the neighbors
$u$ of $v$ in $\Gamma$, equipped 
with the groups $B_u$ 
that generate their product 
$B=\langle B_u\, |\, [u,v] \in 
{\rm edge}(\Gamma)\rangle$ within $G$. 
This product is free because 
there is no edge 
joining two different branches 
$B_u$ and $B_{u'}$, 
hence there is no relation in $B$ 
apart the relations 
defined by the edges of each $B_u$.

Due to Proposition \ref{prop:KS} we have
$G=\varinjlim(\gamma, \cG_\gamma)$,
where $\cG_\gamma$ has 
the vertex groups $A=G_v$
and $B_u$ for each $[u,v] \in 
{\rm edge}(\Gamma)$,
and the common edge group 
$H_v=A\cap B_u$ 
for each edge 
$[u,v] \in {\rm edge}(\gamma)$. 
Moreover, $A$ and $B$
are subgroups of $G$ such that 
$A\cap B=H_v$ and 
$G=A*_{H_v} B$. 
Now the assertion follows from 
Proposition \ref{thm:maximal} applied 
to the latter amalgam. 
\eproof
\subsection{Bearable automorphism groups}
\bdefi[{\rm see~\cite[Definition~2.21]{KPZ17}}] 
\la{def:tree} 
Let $\Gamma$ be a tree and 
$G=\varinjlim(\Gamma, \cG)$ be
a tree product 
with
vertex groups $G_v$, $v \in \vert(\Gamma)$,
where the edge group 
$G_e=G_u\cap G_{v}$ 
is a proper subgroup of $G_v$  
for each
 edge $e=[u,v]$ of $\Gamma$. 
We say that $G$ is \emph{bearable} if
$G_v$ is a nested ind-group for each vertex 
$v\in\vert(\Gamma)$.
\edefi
We have the following analogue 
of Serre's Theorem 
\ref{thm:Serre} for bearable groups.
\bprop[{\rm \cite[Proposition~2.24]{KPZ17}}]
Suppose that the base field $\kk$ is uncountable. 
Let $X$ be an affine variety over $\kk$ 
and $\Gamma$ 
be a tree
with a  countable set $\vert(\Gamma)$. 
Let 
$G=\varinjlim(\Gamma, \cG)$ be 
a bearable subgroup of $\Aut(X)$
such that the vertex groups $G_v$ 
are ind-subgroups of the ind-group $\Aut(X)$.
 Then any algebraic subgroup $H \subset G$ 
 is conjugate to an algebraic subgroup 
 of one of the vertex groups
 $G_v$. 
\eprop
 There are examples of affine surfaces $X$ 
with bearable automorphism 
 groups $\Aut(X)$, 
 see for example  
 \cite{AZ13}, \cite{DG74}--\cite{DG77}, 
 and \cite[Section~2.4]{KPZ17}. 
 In many of these examples, 
 the edge groups $G_e$ 
 are affine algebraic groups
 and the vertex groups $G_v$ 
 are the automorphism groups of  
 $\AA^1$-fibrations.
All automorphism groups considered 
in Sections \ref{sec:A2} and \ref{sec:toric}
are bearable. Indeed, the automorphism 
groups $\Aut(\AA^2)$ and $\Aut_0(\AA^2)$
are  bearable due to the Jung--van 
der Kulk--Nagata theorem.
The groups $\Aut(X_{d,e})$ are bearable due to
Theorem \ref{thm-JvdK-toric-surf}.
%
%
\section{Appendix B: Maximal solvable subgroups 
of $\Aff(\AA^2)$} \la{sec:App-B}
Every solvable 
subgroup of 
an affine algebraic group is contained 
in a maximal solvable 
subgroup, see for example
\cite[Ch.~V, Theorem~10]{Sup76}. 
In this section, we describe
all maximal solvable subgroups of 
$\Aff(\AA^2)$. 
We show that they form exactly 
three distinct conjugacy classes.

We use the following notation. 
Let $Z\subset\GL(2,\kk)$ 
be the subgroup of scalar matrices,
\[\Lambda=\Biggl\{\left(\begin{matrix} 
\lambda & 0\\
0 & \lambda^{-1}\end{matrix}\right), 
\quad \lambda\in\Gm\Biggr\}\subset \SL(2,\kk)
,\quad
\text{and}\quad\tau'=
\left(\begin{matrix} 0 & -1\\
1 & 0\end{matrix}\right)\in \SL(2,\kk).\]
We also let $SB^+\subset\SL(2,\kk)$
be the upper triangular subgroup and  
$2O\subset \SL(2,\kk)$ be the 
binary octahedral group (of order 48). 
\bprop\la{prop:GL2} $\,$
\bnum
\item[$(a)$] The subgroups $SB^+$, 
$\langle \Lambda,\tau'\rangle$, 
and $2O$
are maximal among the 
solvable subgroups
of $\SL(2,\kk)$, and any  maximal
solvable subgroup of $\SL(2,\kk)$
is conjugate to one of them. 

\item[$(b)$] The subgroups $B^+$,
$\langle \TT,\tau\rangle$, and 
$\langle Z, 2O\rangle$ 
are maximal among the 
solvable subgroups
of $\GL(2,\kk)$, and any  maximal
solvable subgroup of $\GL(2,\kk)$
is conjugate to one of them. 

\item[$(c)$] The subgroups 
$\langle B^+, \Trans(\AA^2)\rangle$,
$\langle \TT,\tau, \Trans(\AA^2)\rangle$,
and $\langle Z, 2O, \Trans(\AA^2)\rangle$
are maximal among the 
solvable subgroups
of $\Aff(\AA^2)$, and
 any  maximal
solvable subgroup of $\Aff(\AA^2)$
is conjugate to one of them. 
\enum
\eprop

\bproof Due to the following claim, 
statements (b) and (c) follow from (a).

$\,$

\noindent {\bf Claim 1.} 
\bnum\item[$(a')$]
\emph{If $G\subset\Aff(\AA^2)$ is a 
maximal solvable subgroup,
then 
$G=(Z\times H) \ltimes \Trans(\AA^2)$,
where $H\subset\SL(2,\kk)$
is a  maximal 
solvable subgroup of $\SL(2,\kk)$. 
\item[$(b')$] If $G\subset \GL(2,\kk)$
is a maximal solvable subgroup of $\GL(2,\kk)$, 
then 
$G=Z\times H$, where $H$ is as before. }
\enum

\medskip
\noindent\emph{Proof of Claim 1.} 
The proof is easy and is left to the reader. 
\qed

\medskip

\noindent \emph{Proof of $(a)$.} 
Let us first prove  the second 
assertion of (a). 
Let $H\subset\SL(2,\kk)$ be 
a maximal solvable subgroup.
Then $H$ is closed
in $\SL(2,\kk)$ and is therefore an algebraic
subgroup.  It suffices to consider 
the case where
$H$ is not a Borel subgroup. 
We therefore assume that $H$
is  disconnected. 
Since $H$ is maximal solvable, 
$H$ contains the center 
$C=\{\pm 1\}$ of $\SL(2,\kk)$.

By the Lie-Kolchin theorem, up to conjugation 
in $\SL(2,\kk)$,
we can assume that $H^0$ is contained 
in the upper triangular
subgroup $SB^+\subset \SL(2,\kk)$. 
Since the Borel 
subgroup $SB^+$ is 
maximal solvable in $\SL(2,\kk)$,
we have $H^0\neq SB^+$, 
and therefore
$\dim H^0 \le 1$. 

If $\dim H^0 = 1$, then either 
$H^0\simeq\Gm$ 
is a maximal torus of $SB^+$, 
or $H^0=U^+\simeq\Ga$ 
is the unipotent radical of $SB^+$. 

In the first 
case, we can assume that
$H^0=\Lambda$. 
Being disconnected, $H$ coincides 
with the normalizer 
$\Norm_{\SL(2,\kk)}(\Lambda)
=\langle\Lambda, \tau'\rangle$, 
where $\tau'\notin SB^+$. 
In the second case, $H$ is contained in 
$\Norm_{\SL(2,\kk)}(U^+)=SB^+$. 
Being maximal solvable, $H=SB^+$, 
contrary to our assumption. 

If $H^0=1$, then the finite  subgroup
 $F:=H$ is maximal 
among the solvable finite subgroups of 
$\SL(2,\kk)$. Its image $\bar F$ under
the quotient morphism 
$\rho\colon \SL(2,\kk)\to 
\PGL(2,\kk)$
is maximal among the solvable finite 
subgroups of 
$\PGL(2,\kk)$, and $F=\rho^{-1}(\bar F)$.

Up to conjugacy, $\PGL(2,\kk)$
contains
two infinite series of
solvable  finite subgroups, namely, 
the finite cyclic and  dihedral subgroups. 
Neither of them is maximal solvable. 
There are also exactly 
two conjugacy classes of 
sporadic
solvable finite subgroups 
of $\PGL(2,\kk)$,  
see e.g. 
\cite{Bea10} and \cite{Dol09}.  
The subgroups 
of one are isomorphic to 
the alternating group $A_4$
and of the second to
the symmetric group $S_4$.
Alternating subgroups 
are contained in symmetric subgroups
and are therefore not maximal.
Thus, the only finite subgroups 
that are candidates 
for the maximal solvable finite subgroups of 
$\PGL(2,\kk)$
are isomorphic to $S_4$. 
Their preimages in $SL(2,\kk)$ 
are conjugate to 
the binary octahedral subgroup 
$2O$. 
This proves the second assertion of (a).

To show the second, note that
the Borel subgroup $SB^+$ of $\SL(2,\kk)$
is maximal among the solvable subgroups of 
$\SL(2,\kk)$, cf. Theorem 
\ref{thm:Borel-A2}(e). 
An argument from the proof 
of the second part of (a)
shows that 
$\langle\Lambda, \tau'\rangle$ 
is not conjugate to 
a subgroup of $SB^+$. It is therefore
maximal solvable. 

An argument from the proof above shows
that the binary octahedral subgroup 
$2O$ is maximal among the
solvable  finite subgroups of $\SL(2,\kk)$. 
Since $2O$ is of derived length $4$,
it cannot be conjugate to 
a subgroup of the metabelian
groups $SB^+$ and 
$\langle \Lambda, \tau'\rangle$.
Therefore,
$2O$ is a maximal solvable subgroup 
of $\SL(2,\kk)$.
\eproof
\medskip

{\bf Acknowledgments.} 
We thank St\'ephane Lamy for his patience in 
answering our questions.

\end{document}